%% file: main.tex
\documentclass[12pt, hidelinks]{article}

\usepackage[utf8]{inputenc}
\usepackage{enumerate, enumitem}
\usepackage{url,hyperref}
\usepackage{amsmath,amssymb,amsthm,mathtools}
\usepackage{fullpage}
\usepackage{xcolor}
\usepackage[textsize=small]{todonotes}
\usepackage[small,bf]{caption}
\usepackage{subcaption}
\usepackage{graphicx}
\usepackage{algorithm, algpseudocode}

\graphicspath{{figures/}}

\let\eps\varepsilon

\usepackage[style=alphabetic,backend=bibtex]{biblatex}
\bibliography{refs}

\usepackage[T1]{fontenc}
\usepackage[utf8]{inputenc}

\input defs.tex

\newcommand{\gek}{\succeq}
\todostrue

\title{Dynamic Pricing for Non-fungible Resources:\\
Designing Multidimensional Blockchain Fee Markets}
\author{
    Theo Diamandis\\
    {\small \texttt{tdiamand@mit.edu}}
    \and
    Alex Evans\\
    {\small \texttt{aevans@baincapital.com}}
    \and
    Tarun Chitra\\
    {\small \texttt{tarun@gauntlet.network}}
    \and
    Guillermo Angeris \\
    {\small \texttt{gangeris@baincapital.com}}
}
\date{August 2022}

\newcommand{\fees}{p}

\begin{document}

\maketitle

\begin{abstract}
    Public blockchains implement a fee mechanism to allocate scarce computational 
    resources across competing transactions. 
    Most existing fee market designs utilize a joint, fungible unit of account 
    (\eg, \emph{gas} in Ethereum) to price otherwise non-fungible resources
    such as bandwidth, computation, and storage, by hardcoding their relative prices. 
    Fixing the relative price of each resource in this way inhibits granular price discovery, 
    limiting scalability and opening up the possibility of denial-of-service attacks. 
    As a result, many prominent networks such as Ethereum and Solana have proposed 
    multidimensional fee markets. 
    In this paper, we provide a principled way to design fee markets that 
    efficiently price multiple non-fungible resources. 
    Starting from a loss function specified by the network designer, 
    we show how to compute dynamic prices that align the network’s incentives 
    (to minimize the loss) with those of the users and miners (to maximize their welfare),
    even as demand for these resources changes. 
    Our pricing mechanism follows from a natural decomposition of the network designer’s 
    problem into two parts that are related to each other via the resource prices. 
    These results can be used to efficiently set fees in order to improve network performance.
\end{abstract}

\section{Introduction}

Public blockchains allow any user to submit a \emph{transaction} that modifies the 
shared state of the network. 
Transactions are independently verified and executed by a decentralized network 
of \emph{full nodes}. Because full nodes have finite resources, 
blockchains limit the total computational resources that can be consumed per 
unit of time. 
As user demand may fluctuate, most blockchains implement a \emph{transaction fee} 
mechanism in order to allocate finite computational capacity among competing transactions.

\paragraph{Smart contracts and gas.} 
Many blockchains allow transactions to submit and/or execute programs that
exist on-chain called \emph{smart contracts}. 
Once such a transaction is included in a block, full nodes must re-execute 
the transaction in order to obtain the resulting updated state of the ledger. 
All of these transactions consume computational resources, whose total supply is
finite. To prevent transactions with excessive resource use and transaction spam, 
some smart-contract blockchains require users to pay fees in order to compensate the network for 
processing their transactions.

Most smart contract platforms calculate transaction fees based on a shared unit 
of account. 
In the Ethereum Virtual Machine (EVM), this unit is called \emph{gas}. 
Each operation in the EVM requires a hardcoded amount of gas which is intended 
to reflect its relative resource usage. 
The network enforces a limit on the total amount of gas consumed across all 
transactions in a block. 
This limit, called the \emph{block limit}, prevents the chain from expending 
computational resources too quickly for full nodes to catch up to the latest 
state in a reasonable amount of time. 
Block limits must take into account the maximum amount of each resource that each 
block may consume (such as storage, bandwidth, or memory) without posing an extreme 
burden on full nodes meeting the minimal hardware specifications. 
Because the block limit fixes the total gas supply in each block, the price of
gas in `real' terms (\eg, in terms of US Dollars) fluctuates based on demand
for transactions in the block.

\paragraph{One-dimensional transaction fees.} 
Calculating transaction fees through a single, joint unit of account, such as gas,
introduces two major challenges. 
First, if the hardcoded costs of each operation are not precisely reflective of 
their relative resource usage, there is a possibility of denial-of-service attacks 
(specifically, \emph{resource exhaustion attacks}~\cite{perez2019broken}), where 
an attacker exploits resource mispricing to overload the network.  
Historically, the Ethereum network has suffered from multiple DoS
attacks~\cite{suicide-attack, transaction-attack, wilcke_ethereum_2016} and has
had to manually adjust the relative prices accordingly
(\eg,~\cite{eip150,eip2929}). Amending the hardcoded costs of each gas
operation in response to such attacks typically requires a hard fork of the
client software. 

Second, one-dimensional fee markets limit the theoretical throughput of the network. 
Using a joint unit of account to price separate resources decouples their price 
from supply and demand. As an extreme example to illustrate this dynamic, 
if the block gas limit is fully saturated with exclusively CPU-intensive operations, 
gas price will increase as transactions compete for limited space. 
The cost of transactions that consume exclusively network bandwidth 
(and nearly no CPU resources) will also increase because these also require gas, 
even if demand and supply for bandwidth resources across the network remain unchanged. 
As a result, fewer bandwidth-intensive transactions can be included in the block 
despite spare capacity, limiting throughput. 
This limitation occurs because the shared unit of account only allows the network 
to price resources \emph{relative to each other} and not in real terms based on 
the supply and demand for each resource. 
As we will soon discuss, allowing resources to be priced separately promotes more 
efficient resource utilization by enabling more precise price discovery. 
We note that this increase in throughput need not increase hardware requirements 
for full nodes.

\paragraph{Multidimensional fee markets.} 
Due to the potential scalability benefits of more granular price discovery, 
a number of popular smart contract platforms are actively exploring multidimensional 
fee market mechanisms~\cite{adler2022ethcc,multidimensional-eip1559}. 
We discuss some example proposals that are under active development, below.

\paragraph{Rollups and data markets.} 
Rollups are a popular scaling technology that effectively decouples transaction 
validation and execution from data and consensus~\cite{buterin_incomplete_2021}. 
In rollups, raw transaction data is posted to a base blockchain. 
Rollups also periodically post succinct proofs of valid execution to the base 
chain in order to enable secure bridging, prevent rollbacks, and arbitrate potential 
disputes by using the base chain as an anchor of trust. 
Rollups naturally create two separate fee markets, one for base layer transactions 
and one for rollup execution. As rollups have become a popular design pattern 
for achieving scalability, specialized blockchains (called \emph{lazy blockchains}) 
that exclusively order raw data through consensus (\ie, do not perform execution) 
have emerged~\cite{al-bassam_fraud_2018,al-bassam_lazyledger_2019,
polygon_team_introducing_2021, nazirkhanova_information_2021, tas_light_2022}. 
These blockchains naturally allow for transaction data/bandwidth and execution to 
be priced through independent (usually one-dimensional) fee markets~\cite{adler2021ethcc}. 
Similarly, Ethereum developers have proposed EIP-2242, wherein users may submit 
special transactions which contain an additional piece of data called a 
\emph{blob}~\cite{adler2019dataavail,adler2019eip2242}.
Blobs may contain arbitrary data intended to be interpreted and executed by 
rollups rather than the base chain.
Later, EIP-4844 extended these ideas by establishing a two-dimensional fee market 
wherein data blobs and base-chain gas have different limits and are priced 
separately~\cite{noauthor_proto-danksharding_2022}. 
EIP-4844 therefore intends to increase scalability for rollups, as blobs do not 
have to compete with base-chain execution for gas.

\paragraph{Incentivizing parallelization.} 
Most smart contract platforms, including the EVM, execute program operations 
sequentially by default, limiting performance.
There are several proposals to enable parallel execution in the EVM which 
generally fall into two categories. 
The first involves minimal changes to the EVM and pushes the responsibility of 
identifying opportunities for parallel execution to full nodes~\cite{
    gelashvili_block-stm_2022,chen_forerunner_2021,saraph_empirical_2019}.
The other approach involves \emph{access lists} which require users to specify 
which accounts their transaction will interact with, allowing the network to 
easily identify non-conflicting transactions that can be executed in 
parallel~\cite{vbuterin_easy_2017, buterin_eip-2930_2020}.
While Ethereum  makes access lists optional, other virtual machine implementations, 
such as Solana Sealevel and FuelVM, require users to specify the accounts their 
transaction will interact with~\cite{yakovenko_sealevel_2020, 
bvrooman_github_nodate, adler_accounts_2020}.
Despite this capability, a large fraction of transactions often want to access 
the same accounts in scenarios including auctions, arbitrage 
opportunities, and new product launches.
Such contention significantly limits the potential benefits of the virtual 
machine's parallelization capabilities.
As a result, developers of Solana have proposed multidimensional fee markets 
that price interactions with each account separately in order to charge higher 
fees for transactions which require sequential execution~\cite{aeyakovenko_consider_2021}. 
Such a proposal incentivizes usage of spare capacity on multi-core machines.

\paragraph{This paper.} 
In this paper, we formally illustrate how to efficiently update resource prices, 
what optimization problem these updates attempt to solve, 
and some consequences of these observations. 
We also numerically demonstrate that this approach enhances network performance 
and reduces DoS-style resource exhaustion attacks. 
We frame the pricing problem in terms of an idealized, omniscient network designer 
who chooses transactions to include in blocks in order to maximize total welfare, 
subject to demand constraints. 
(The designer is omniscient as the welfare is unknown and likely unmeasurable 
in any practical setting.) 
We show that this problem, which is the ‘ideal end state’ of a blockchain but 
not immediately useful by itself, decomposes into two problems, 
coupled by the resource prices. 
One of these two problems is a simple one which can be easily solved on chain 
and represents the cost to the network for providing certain resources, 
while the other is a maximal-utility problem that miners and users implicitly 
solve when creating and including transactions for a given block. 
Correctly setting the resource prices aligns incentives such that the resource 
costs to the network are exactly balanced by the utility gained by the users and miners. 
This, in turn, leads to block allocations which solve the original ‘ideal’ problem, on average.

For convenience, we provide appendix~\ref{app:convex} as a short introduction
to convex optimization. We recommend readers unfamiliar with convex
optimization at least skim this appendix, as it provides a short introduction
to all the mathematical definitions and major theorems used in this paper. As a
general guideline, we recommend those uninterested in theoretical results to
skip~\S\ref{sec:min-demand} and~\S\ref{sec:properties} on a first reading.

\subsection{Related work}
The resource allocation problem has been studied in many fields, including
operations research and computer systems. Agrawal, et
al.~\cite{agrawal2022allocation} proposed a similar formulation and price
update scheme for fungible resources where utility is defined per-transaction.
Prior work on blockchain transaction fees varies from the formal axiomatic
analysis of game theoretic properties that different fee markets should
have~\cite{roughgarden2020eip1559, chung2021foundations}
to analysis of dynamic fees from a direct algorithmic
perspective~\cite{ferreira2021dynamic, leonardos2021dynamical,
reijsbergen2021transaction}. Works of the latter form generally focus on
whether the system macroscopically converges to an equilibrium. Moreover, these
mechanisms focus on dynamic pricing at the block level (\eg, how many
transactions should be allowed in a block?) and not directly on questions of
how capacity should be allocated and priced across different transaction types.

\paragraph{EIP-1559.}
EIP-1559~\cite{eip1559}, implemented last year, proposed major changes to
Ethereum's transaction fee mechanism. Specifically, EIP-1559 implemented a base
fee for transactions to be included in each block, which is dynamically
adjusted to hit a target block usage and burned instead of rewarded to the
miners. We note that while EIP-1559 is closely related to the problem we
consider, it ultimately has a different goal: EIP-1559 attempts to make the fee
estimation problem easier in a way that disincentivizes manipulation and
collusion~\cite{vitalikfirstpricedauctions,roughgarden2020eip1559}. We, on the
other hand, aim to price resources dynamically to achieve a given
network-specified objective. Finally, we note that prior work such
as~\cite{ferreira2021dynamic} has proved incentive compatibility for a large
set of mechanisms that are a superset of EIP-1559. It is likely (but not proven
in this work) that our model fits within an extension of their incentive
compatibility framework. We leave game theoretic analysis and strategies to
ensure incentive compatibility as future work.

\section{Transactions and resources}\label{sec:tx-model}
Before introducing the network's resource pricing problem, we discuss the
general set up and motivation for the problem in the case of blockchains.

\paragraph{Transactions.} We will start by reasoning about \emph{transactions}.
A transaction can represent arbitrary data sent over the peer-to-peer network
in order to be appended to the chain. Typically, a transaction will represent a
value transfer or a call to a smart contract that exists on chain. These
transactions are broadcasted by users through the peer-to-peer network and
collected by consensus nodes in the \emph{mempool}, which contains all submitted
transactions that have not been included on chain. A \emph{miner} gets to choose
which transactions from the mempool are included on chain. Miners may also
outsource this process to a \emph{block builder} in exchange for a
reward~\cite{buterinPBS}. Once a transaction is included on chain, it is
removed from the mempool.
(Any conflicting transactions are also removed from the mempool.)

\paragraph{Nodes.} Every transaction needs to be executed by \emph{full nodes}
(which we will refer to as `nodes'). Nodes compute the current state of the chain
by executing and checking the validity of all transactions that have been included on the chain
since the genesis block.
Many blockchains have minimum computational
requirements for nodes in a blockchain: any node meeting these requirements
should be able to eventually `catch up' to the head of the chain in a reasonable
amount of time, \ie, execute all transactions and reach the latest state, even
as the chain continues to grow. 
(For example, Ethereum requires 4GB RAM and 2TB of SSD Storage, and
they recommend at least a Intel NUC, 7th gen processor~\cite{ethnodereqs}.)
These requirements both limit the computational resources each block is 
allowed to consume and promote decentralization by ensuring the required hardware
does not become prohibitively expensive.
If transactions are included in a blockchain faster than nodes are able to execute
them, nodes cannot reconstruct the latest state and can't ascertain 
the validity of the chain. 
This type of denial of service (DoS) attack is also referred 
to as a resource exhaustion attack.
(As a side note, in some systems, it is possible to provide an
easily-verifiable proof that the state is correct without the need to execute
all past transactions to validate the state of the chain. In these systems, the
time-consuming step is generating the proofs. A similar market mechanism might
make sense for these systems, but we do not explore this topic here.)

\paragraph{Resource targets and limits.} There are several ways to prevent this
type of denial of service attack. For example, one method is to enforce that
any valid transaction (or sequence of transactions, \eg, a block) consumes at
most some fixed upper bound of resources, or combinations of resources.
These limits are set so
that that a node satisfying the minimum node requirements is always able to
process transactions quickly enough to catch up to the head of the chain in a
reasonable amount of time. Another possibility is to disincentivize miners and
users from repeatedly including transactions that consume large amounts of
resources while allowing short `bursts' of resource-heavy transactions. This
margin needs to be carefully balanced so that a node meeting the minimum
requirements is able to catch up after a certain period of time. This intuition
suggests having both a  `resource target' and a larger `resource limit,' which
we will formalize in what follows.

\paragraph{Resources.} Most blockchain implementations define a
number of \emph{meterable resources} (simply \emph{resources} from here on out)
which we will label $i=1, \dots, m$, that a transaction can consume. For
example, in Ethereum, the resources could be the individual Ethereum
Virtual Machine (EVM) opcodes used in the execution of a transaction. In this
paper, the notion of a `resource' is much more general than simply an `opcode'
or an `execution unit'. Here, a resource can refer to anything as coarse as
`total bandwidth usage' or as granular as individual instructions executed on a
specific core of a machine---the only requirement for a resource, as used in
this paper, is that it can be easily and consistently metered across any node. For
a given transaction $j=1, \dots, n$, we will let $a_j \in \reals^m_+$ be the
vector of resources that transaction $j$ consumes. In particular, the $i$th
entry of this vector, $(a_j)_i$, denotes the amount of resource $i$ that
transaction $j$ uses.
We note that the resource usage $(a_j)_i$ does not, in fact,
need to be nonnegative.
While our mechanism works in the more general case (with some small modifications),
we assume nonnegativity in this work for simplicity.

\paragraph{Combined resources.} This framework naturally includes combinations 
of resources as well.
For example, we may have two resources $R_1$ and $R_2$, each
cheap to execute once in a transaction, which are costly to execute serially 
(\ie, it is costly to execute $R_1$ and then $R_2$ in the same transaction). 
In this case, we
can create a `combined' resource $R_1R_2$ which is itself metered separately.
Note that, in our discussion of resources, there is no requirement that the
resources be independent in any sense and such `combined resources' are
themselves very reasonable to consider.

\paragraph{Resource utilization targets.} As mentioned previously, many
networks have a minimum node requirement, implying a sustained target for
resource utilization in each group of transactions added to the blockchain. 
(For simplicity, we will call this sequence of
transactions a \emph{block}, though it could be any collection of transactions
that makes sense for a given blockchain.) We will write this resource
utilization target as a vector $b^\star \in \reals^m$ whose $i$th entry denotes
the desired amount of resource $i$ to be consumed in a block. The resource
utilization of a particular block is a linear function of the transactions
included in a block, written as a Boolean vector $x \in \{0, 1\}^n$, whose
$j$th entry is one if transaction $j$ is included in the block and is zero
otherwise. We will write $A \in \reals^{m\times n}$ as a matrix whose $j$th
column is the vector of resources $a_j$ consumed by transaction $j$. We can
then write the total quantity of consumed resources by this block as
\begin{equation}\label{eq:quantity}
    y = Ax,
\end{equation}
where $y \in \reals^m$ is a vector whose $i$th entry denotes the quantity of
resource $i$ consumed by all transactions included in the block. Additionally,
we can write the deviation from the target utilization, sometimes called the
residual, as
\[
    Ax - b^\star,
\]
\ie, a vector whose $i^\mathrm{th}$ element is the total quantity of resource
$i$, consumed by transactions in this block, minus the target for resource $i$,
$b_i^\star$. For example, in Ethereum post EIP-1559, there is only one resource,
gas, which has a target of 15M~\cite{ethgas}. (We will see later how this notion of a
`resource utilization target' can be generalized to a loss function.)

\paragraph{Resource utilization limits.}
In addition to resource targets, a blockchain may introduce a resource limit $b$
such that any valid block with transactions $x$ must satisfy
\[
    Ax \leq b.
\]
Continuing the example from before, Ethereum after EIP-1559 has a single resource, gas,
with a resource limit of 30M.

\paragraph{Network fees.} 
We want to incentivize users and miners to keep the total resource usage
`close' to $b^\star$. To this end, we introduce a \emph{network fee} $\fees_i$
for resource $i = 1,\dots, m$, which we will sometimes call the
\emph{resource price}, or just the \emph{price}. If transaction $j$ with
resource vector $a_j$ is included in a block, a fee $\fees^Ta_j = \sum_i
\fees_i (a_j)_i$ is paid to the network. (In Ethereum, the network fee is
implemented by burning some amount of the gas fee for each transaction and can
be thought of as the \emph{burn rate}.) For now, we will assume that as the fee
$\fees_i$ gets larger, the amount of resource $i$ used in a block $(Ax)_i$ will
become smaller and vice versa.

\paragraph{Resource mispricing.}
Given $A$ and $b$, it is, in general, not clear how to set the fees $\fees$ in
order to ensure the network has good performance. (In other words, so that the
resource utilization is `close' to $b^\star$.) As a real world example,
starting in Ethereum block number $2283416$ (Sept.\ 18, 2016) an attacker
exploited the fact that resources were mispriced for the EXTCODESIZE opcode,
causing the network to slow down meaningfully. This mispricing was fixed via
the hard fork on block $2463000$ (Oct.\ 18, 2016) with details outlined in
EIP-150~\cite{eip150}. (The effects of this mispricing can still be observed
when attempting to synchronize a full node. A dramatic slowdown in downloading
and processing blocks happens starting at the previously mentioned block
height.) Though usually less drastic, there have been similar events on other
blockchains, underscoring the importance of correctly setting resource prices. 

\paragraph{Setting fees.}
We want a simple update rule for the network fees $\fees$ with the following properties:
\begin{enumerate}
    \item If $Ax = b^\star$, then there is no update.
    \item If $(Ax)_i > b_i^\star$, then $\fees_i$ increases.
    \item Otherwise, if $(Ax)_i < b_i^\star$, then $\fees_i$ decreases.
\end{enumerate}
There are many update rules with these properties. As a simple example, we can
update the network fees as
\begin{equation}\label{eq:basic-update}
    \fees^{k+1} = \left(\fees^k + \eta(Ax - b^\star)\right)_+,
\end{equation}
where $\eta > 0$ some (usually small) positive number, often referred to as the
`step size' or `learning rate', $k$ is the block number such that $p^k$ are
the resource prices at block $k$, and
$(w)_+ = \max\{0, w\}$ for scalar $w$ and is applied elementwise for vectors.
Recently, Ethereum developers~\cite{vitalikmultidimeip1559} proposed the update rule
\begin{equation}\label{eq:buterinUpdate}
    \fees^{k+1} = \fees^k \odot \exp\left(\eta(Ax - b^\star)\right).
\end{equation}
Here, $\exp(\cdot)$ is understood to apply elementwise, while $\odot$ is the
elementwise or Hadamard product between two vectors. The remainder of this
paper will show that many update rules of this form are attempting to
(approximately) solve an instance of a particular optimization problem with a
natural interpretation, where parts of the update rule come from a specific
choice of objective function by the network designer. (We show in
appendix~\ref{app:exp-update} that a similar rule to~\eqref{eq:buterinUpdate}
can be derived as a consequence of the framework presented in this paper, under
a particular choice of variable transformation.)

We note that~\cite{ferreira2021dynamic} has studied fixed points of such
iterations in the one-dimensional case, extending the analysis of EIP-1559.
However, the multidimensional scenario can be quite a bit more subtle to
analyze. For instance, the multiplicative update rule~\eqref{eq:buterinUpdate}
can admit `vanishing gradient' behavior in high-dimensions
\cite{hochreiter1998vanishing}. We suspect that the one-dimensional fee model
of~\cite{ferreira2021dynamic} can be extended to the multidimensional
rules~\eqref{eq:basic-update} and~\eqref{eq:buterinUpdate} and leave this for
future work.

\section{Resource allocation problem}
As system designers, our ultimate goal is to maximize utility of the underlying
blockchain network by appropriately allocating resources to transactions.
However, we cannot perform this allocation directly, since we cannot control
what users nor miners wish to include in blocks, nor do we know what the
utility of a transaction is to users and miners. Instead, we aim to set the
fees $\fees$ to ensure that the resource usage is approximately equal to the
desired target, which we will represent as an objective function. We will show
later that the update mechanisms proposed above naturally fall out of a more
general optimization formulation. Similarly to the Transmission Control Protocol (TCP),
each update rule is a
result of a particular objective function, chosen by the network 
designer~\cite{low1999optimization, low2003duality}.

\paragraph{Loss function.}
We define a \emph{loss function} of the resource usage, $\ell: \reals^m \to
\reals \cup \{\infty\}$, which maps a block's resource utilization, $y$, to the
`unhappiness' of the network designer, $\ell(y)$. We assume only that $\ell$ is
convex and lower semicontinuous. (We will not require monotonicity,
nonnegativity, or other assumptions on $\ell$, but we will show that useful
properties do hold in these scenarios.)

For example, the loss function can encode `infinite dissatisfaction' if the
resource target is violated at all:
\begin{equation}\label{eq:indicator-loss}
    \ell(y) = \begin{cases}
        0 & y = b^\star \\
        \infty & \text{otherwise}.
    \end{cases}
\end{equation}
(Functions of this form, which are either 0 or $\infty$ at every point, are
known as \emph{indicator functions}.) Note also that this loss is not
differentiable anywhere, but it is convex. Another possible loss, which is also an
indicator function, is
\begin{equation}\label{eq:inequality-loss}
\ell(y) = \begin{cases}
    0 & y \le b^\star\\
    \infty & \text{otherwise}.
\end{cases}
\end{equation}
This loss can roughly be interpreted as: we do not mind any usage below
$b^\star$, but we are infinitely unhappy with any usage above the target amounts.
Alternatively, we may only care about large deviations from the target
$b^\star$:
\[
    \ell(y) = (1/2)\|y - b^\star\|_2^2,
\]
or, perhaps, require that the loss is simply linear and independent of $b^\star$,
\begin{equation}\label{eq:linear-loss}
    \ell(y) = u^Ty,
\end{equation}
for some fixed vector $u \in \reals^m$. Another important family of losses are
those which are separable and depend only on the individual resource utilizations,
\begin{equation}\label{eq:nondec-loss}
\ell(y) = \sum_{i=1}^m \phi_i(y_i)
\end{equation}
where $\phi_i: \reals \to \reals \cup \{\infty\}$ for $i=1, \dots, m$, are
convex, nondecreasing functions. (The loss~\eqref{eq:inequality-loss} is a
special case of this loss, while~\eqref{eq:linear-loss} is a special case when
the vector $u$ is nonnegative.) We will make the technical assumption that
$\phi_i(0) < \infty$ for every $i$, otherwise no resource allocation would have
finite loss.

We will see that each definition of a loss function implies a particular update rule
for the network fees $\fees$.
This utility function can more generally be engineered to appropriately capture
tradeoffs in increasing throughput of a particular resource at the possible
detriment of other resources.

\paragraph{Resource constraints.} Now that we have defined the network
designer's loss, which is a way of quantifying `unhappiness' when the resource
usage is $y$, we need some way to define the transactions that users are willing to
create and, conversely, that miners are willing (and able) to include. We do this
in a very general way by letting $S \subseteq \{0, 1\}^n$ be the set of
possible transactions that users and miners are willing and able to create and include.
Note that this set is discrete and can be very complicated or potentially
hard to maximize over (as is the case in practice). For example, the set $S$
could encode a demand for transactions which depend on other transactions being
included in the block (as is the case in, \eg, miner extractable
value~\cite{kulkarniMEV, daian2020flash, qin2021quantifying}), network-defined
hard constraints of certain resources (such as $Ax \le b$ for every $x \in S$),
and even very complicated interactions among different transactions
(if certain contracts can, for example, only be called a fixed number of times,
as in NFT mints). We make no assumptions about the structure of this set $S$,
but only require that the included transactions, $x \in \{0, 1\}^n$, obey the
constraint $x \in S$.

\paragraph{Convex hull of resource constraints.} 
A network designer may be more interested in the long-term resource utilization of the network
than the resource utilization of any one particular block.
In this case, the designer may choose to `average out' each transaction over a number of
blocks instead of deciding whether or not to include it in the next block.
To that end, we,
as designers, will be allowed to choose convex combinations of elements of the
constraint set $S$, which we will write as $\conv(S)$. (In general, this means
that we can pick probability distributions over the elements of $S$, and $x$ is
allowed to be the expectation of any such probability distribution; \ie, we
only require that, for the designer, $x$ is reasonable `in expectation'.) 
Here, components of $x$ may vary continuously between $0$ and $1$; these
values have a simple interpretation. If $x_i$ is not $0$ or $1$, then we can
interpret the quantity $1/x_i$ as only including transaction $i$ after
roughly $1/x_i$ blocks. Of course, neither users nor miners can
choose transactions to be `partially included', so this property will only
apply to the idealized problem we present below. While this relaxation might
seem unrealistically `loose', we will see later how this easily translates to
the realistic case where transactions are either included or not (that is,
$x_i$ is either $0$ or $1$) by users and miners.

\paragraph{Transaction utility.} Finally, we introduce
the \emph{transaction utilities}, which we will write as $q \in \reals^n$.
The transaction utility $q_j$ for transaction $j=1, \dots, n$ denotes the users'
and miners' joint utility of including transaction $j$ in a given block. Note that
it is very rare (if at all possible) to know the values of $q$. However, we will see that,
under mild assumptions, we do not need to know the values of $q$ in order
to get reasonable prices, and reasonable update rules will not depend on $q$.

\paragraph{Resource allocation problem.}
We are now ready to write the \emph{resource allocation problem}, which 
is to maximize the utility of the included transactions, minus the loss,
over the convex hull of possible transactions:
\begin{equation}\label{eq:allocation}
    \begin{aligned}
        & \text{maximize}     && q^Tx - \ell(y) \\
        & \text{subject to}   && y = Ax \\
                            &&& x \in \conv(S).
    \end{aligned}
\end{equation}
This problem has variables $x \in \reals^n$ and $y \in \reals^m$, and the
problem data are the resource matrix $A \in \reals^{m\times n}$, the set of
possible transactions $S\subseteq \{0, 1\}^n$, and the transaction utilities $q
\in \reals^n$. Because the objective function is concave and the constraints
are all convex, this is a convex optimization problem (see appendix~\ref{app:convex}).
On the other hand, even
though the set $\conv(S)$ is convex, it is possible that $\conv(S)$ does not
admit an efficient representation (for example, it may contain exponentially
many constraints) which means that solving this problem is, in general, not easy.

\paragraph{Interpretation.} We can interpret the resource allocation
problem~\eqref{eq:allocation} as the `best case scenario', where the designer is
able to choose which transactions are included (or `partially included') in a
block in order to maximize the total utility. While this problem is not
terribly useful by itself, since (a) it cannot really be implemented in
practice, (b) we often don't know $q$, and (c) we cannot `partially include' a
transaction within a block, we will see that it will decompose naturally into
two problems. One of these problems can be easily solved on chain, while the
other is solved implicitly by the users (who send transactions to be included)
and miners (who choose which transactions to include). The solutions to the latter
problem can always be assumed to be integral; \ie, no transactions are
`partially included'. This will allow us to construct a simple update rule for
the prices, which does not depend on $q$. For the remainder of the paper, we
will call this combination of users and miners the \emph{transaction
producers}.

\paragraph{Offchain agreements and producers.} Due to the inevitability of
user-miner collusion, we consider the combination of the two, the transaction
producers, as the natural unit. For example, it is not easily possible to
create a transaction mechanism where the users are forced to pay miners some
fixed amount, since it is always possible for miners to refund users via some
off-chain agreement~\cite{roughgarden2020eip1559}. Similarly, we cannot force
miners to accept certain transactions by users, since a miner always has
plausible deniability of not having seen a given transaction in the mempool.
While not perfect for a general analysis of incentives, this conglomerate
captures the dynamics between the network's incentives and those of the miners
and users better than assuming each is purely selfishly maximizing their own
utility (as opposed to strategically colluding) and suffices for our purposes.

\subsection{Setting prices using duality}
In this section, we will show a decomposition method for this problem. This
decomposition method suggests an algorithm (presented later) for iteratively
updating fees in order to maximize the transaction producers' utility minus the
loss of the network, given historical observations.

To start,
we will reformulate~\eqref{eq:allocation} slightly by pulling the constraint $x
\in \conv(S)$ into the objective,
\begin{equation}\label{eq:allocation-ind}
    \begin{aligned}
        & \text{maximize}     && q^Tx - \ell(y) - I(x) \\
        & \text{subject to}   && y = Ax
    \end{aligned}
\end{equation}
where $I: \reals^n \to \reals \cup \{\infty\}$ is the indicator function defined as
\[
I(x) = \begin{cases}
    0 & x \in \conv(S)\\
    +\infty & \text{otherwise}.
\end{cases}
\]

\paragraph{Dual function.}
The Lagrangian~\cite[\S 5.1.1]{boyd2004convex} for problem~\eqref{eq:allocation-ind} is then
\[
    L(x, y, \fees) = q^Tx - \ell(y) - I(x) + \fees^T(y - Ax),
\]
with dual variable $\fees \in \reals^m$. This corresponds to `relaxing' the 
constraint $y = Ax$ to a penalty $p^T(y - Ax)$ assigned to the objective, where
the price per unit violation of constraint $i$ is $p_i$. (Negative values denote
refunds.) Rearranging slightly, we can write
\[
    L(x, y, \fees) = \fees^Ty -\ell(y) + (q - A^T\fees)^Tx - I(x).
\]
The corresponding dual function~\cite[\S 5.1.2]{boyd2004convex},
which we will write as $g: \reals^m \to \reals \cup \{-\infty\}$, is found by maximizing
over $x$ and $y$:
\begin{equation}\label{eq:dual-fn}
g(\fees) = \sup_y \left(\fees^Ty - \ell(y)\right) + \sup_x \left((q - A^T\fees)^Tx - I(x)\right).
\end{equation}

\paragraph{Discussion.}
The first term can be recognized as the Fenchel conjugate of
$\ell$~\cite[\S 3.3]{boyd2004convex} evaluated at $p$, which we will write as $\ell^*(p)$,
while the second term is the optimal value of the following problem:
\begin{equation}\label{eq:miner-packing}
\begin{aligned}
    & \text{maximize} && (q - A^T\fees)^Tx\\
    & \text{subject to} &&x \in \conv(S),
\end{aligned}
\end{equation}
with variable $x \in \reals^n$. We can interpret this problem as the
transaction producers' problem of creating and choosing transactions to be
included in a block in order to maximize their utility, after netting the fees
paid to the network. We note that the optimal value of~\eqref{eq:miner-packing}
in terms of $\fees$, which we will write as $f(p)$, is the pointwise supremum
of a family of linear (and therefore convex) functions of $\fees$, so it, too,
is a convex function~\cite[\S 3.2.3]{boyd2004convex}. Note that since the
objective is linear, problem~\eqref{eq:miner-packing} has the same optimal
objective value as the nonconvex problem
\[
    \begin{aligned}
        & \text{maximize} && (q - A^T\fees)^Tx\\
        & \text{subject to} &&x \in S,
    \end{aligned}
\]
where we have replaced $\conv(S)$ with $S$. (See appendix~\ref{app:convex-hull}
for a simple proof.)  Finally, since $f(\fees)$ is the optimal value of
problem~\eqref{eq:miner-packing} for fees $\fees$, the dual function $g$ can be
written as
\[
g(\fees) = \underbrace{\ell^*(\fees)}_{\text{network}} + 
\underbrace{f(\fees)}_{\text{tx producers}}.
\]
Since the dual function $g$ is the sum of convex functions $\ell^*$ and $f$,
it, too, is convex. (We will make use of this property soon.) Having defined
the dual function $g$, we will see how this function can give us a criterion
for how to best set the network fees $\fees$.

\paragraph{Duality.} An important consequence of the definition of the dual
function $g$ is \emph{weak duality}~\cite[\S 5.2.2]{boyd2004convex}. Specifically,
letting $s^\star$ be the optimal objective value for
problem~\eqref{eq:allocation}, we have that
\[
    g(p) \ge s^\star,
\]
for every possible choice of price $p\in \reals^m$. This is true because we have
essentially `relaxed' the constraint to a penalty, so any feasible point $x, y$
for the original problem~\eqref{eq:allocation-ind} always has $0$ penalty.
(There may, of course, be other points that are not feasible
for~\eqref{eq:allocation-ind} but are perfectly reasonable for this `relaxed'
version, so we've only made the set of possibilities larger.) The proof is a
single line:
\[
    g(p) = \sup_{x, y} L(x, y, p) \ge \sup_{y = Ax} L(x, y, p) 
    = \sup_{y =Ax} \left(q^Tx - \ell(y) - I(x)\right) = s^\star.
\]
A deep and important result in convex optimization is that, in fact, there
exists a $\fees^\star$ for which
\[
    g(p^\star) = s^\star,
\]
under some basic constraint qualifications.\footnote{The condition is that the
    relative interior of $A\conv(S) \cap \dom \ell$ is nonempty. Here, we write
    $A\conv(S) = \{Ax \mid x \in \conv(S)\}$ and $\dom \ell = \{x \mid \ell(x)
    < \infty\}$, while the relative interior is taken with respect to the
affine hull of the set. This condition almost always holds in practice for
reasonable functions $\ell$ and sets $S$.} In other words, adding the
constraint $y = Ax$ to the problem is identical to correctly setting the prices
$p$. Since we know for any $p$ that $g(p) \ge s^\star$ then
\[
    g(p^\star) = \inf_p g(p),
\]
or, that $p^\star$ is a minimizer of $g$. This motivates an optimization
problem for finding the prices.

\paragraph{The dual problem.}
The \emph{dual problem} is to minimize $g$, as a function of the fees $\fees$.
In other words, the dual problem is to find the optimal value of
\begin{equation}\label{eq:dual-problem}
\begin{aligned}
    & \text{minimize} && g(\fees),
\end{aligned}
\end{equation}
with variable $\fees \in \reals^m$. If we can easily evaluate $g$, then, since
this problem is a convex optimization problem, as $g$ is convex, solving it is
usually also easy. An optimizer of the dual problem has a simple interpretation
using its optimality conditions. Let $\fees^\star$ be a solution to the dual
problem~\eqref{eq:dual-problem} for what follows. If the packing
problem~\eqref{eq:miner-packing} has a unique solution $x^\star$ for
$\fees^\star$, then the objective value $f$ is differentiable at $\fees^\star$.
(See appendix~\ref{app:convex-hull}.) Similarly, under mild conditions on the
loss function $\ell$ (such as strict convexity) the function $\ell^*$ is
differentiable at $\fees^\star$, with derivative $y^\star$ satisfying $\nabla
\ell(y^\star) = p^\star$. In this case, the optimality conditions for
problem~\eqref{eq:dual-problem} are that
\begin{equation}\label{eq:grad}
    \nabla g(\fees^\star) = y^\star - Ax^\star = 0.
\end{equation}
In other words, the fees $\fees^\star$ that minimize~\eqref{eq:dual-problem}
are those that charge the transaction producers the exact marginal costs faced
by the network, $\nabla \ell(Ax^\star) = \fees^\star$. Furthermore, these are
exactly the fees which incentivize transaction producers to include transactions that maximize
the welfare generated minus the network loss, subject to
resource constraints, since $y^\star$ and $x^\star$ are feasible and optimal
for problem~\eqref{eq:allocation}.

\paragraph{Differentiability.} In general, $g$ is not always differentiable,
but is almost universally subdifferentiable, under mild additional conditions
on $\ell$ (\eg, $\ell$ does not contain a line). Condition~\eqref{eq:grad} may
then be replaced with
\[
    0 \in \partial g(p^\star) = -Y^\star(p^\star) + AX^\star(p^\star),
\]
where
\[
    Y^\star(p) = \argmax_y \left(p^Ty - \ell(y)\right),
\]
while $X^\star(p)\subseteq \conv(S)$ are the maximizers of
problem~\eqref{eq:miner-packing} for price $p$. We define $AX^\star(p) = \{Ax
\mid x \in X^\star(p)\}$, and write $\partial g(p^\star)$ for the subgradients
of $g$ at $p^\star$ (\cf, appendix~\ref{app:convex-subdiff}). The condition then says that the
intersection of the extremizing sets $Y^\star(p^\star)$ and $AX^\star(p^\star)$
is nonempty at the optimal prices $p^\star$. We show a special case of this
below, when $p=0$, with a direct proof using strong duality that does not
require these additional conditions.

\subsection{Minimal demand conditions}\label{sec:min-demand}
We can give a condition for which we can guarantee that the optimal prices,
\ie, those which minimize the dual problem~\eqref{eq:dual-problem}, satisfy $p \ne
0$. The condition is: when resources have zero fee, the optimal
set of included transactions that would be included at no price, defined as
$X^\star \subseteq [0, 1]^n$, with
\[
    X^\star = \argmax_{x \in \conv(S)} \,q^Tx,
\]
is `disjoint' from the set of minimizers of the loss, $Y^\star \subseteq \reals^m_+$,
defined
\[
    Y^\star = \argmin_{y} \,\ell(y),
\]
in the following sense:
\[
    AX^\star \cap Y^\star = \emptyset,
\]
where $AX^\star = \{Ax \mid x \in X^\star\}$. An intuitive version of this
condition is that the demand for transactions, if they could be executed at no
cost to the transaction producers, always incurs some loss for the network.
This, in turn, implies that the optimal fees $p$ for such resources must be
nonzero.

For convenience, for the rest of this section, we will define $f^*(x) = q^Tx$
when $x \in \conv(S)$ and $-\infty$ otherwise. This lets us write:
\[
    X^\star = \argmax f^*,
\]
and, for any $x^\star \in X^\star$, we have
\[
    \sup_{x \in \conv(S)} q^Tx = q^Tx^\star = f^*(x^\star) = \sup f^*.
\]

\paragraph{Proof.} To see this, we will use strong duality. We will prove the
contrapositive statement: if $g(0)$ is optimal, then there exists a point in
the intersection $AX^\star \cap Y^\star$.

If $g(0)$ is optimal, then, using
strong duality:
\[
    g(0) = \sup_{\substack{x \in \conv(S)\\ y=Ax}} q^Tx - \ell(y) = \sup_{y = Ax} f^*(x) - \ell(y).
\]
Rewriting the problem to remove the constraint $y=Ax$, we have
\[
    g(0) = \sup_x f^*(x) - \ell(Ax)
\]
Since $\conv(S)$ is a compact set and $\ell$ is lower semi-continuous, there
exists some $\tilde{x} \in \conv(S)$ which achieves this maximum. Now, using the
definition of $g$ (\cf, equation~\eqref{eq:dual-fn}),
\[
    g(0) = \sup_{x, y} f^*(x) - \ell(y) = \sup f^* - \inf \ell.
\]
Putting both statements together, we find
\[
    f^*(\tilde{x}) - \ell(A\tilde{x}) = \sup f^* - \inf \ell.
\]
Since we know, by definition of $\sup$ and $\inf$ that
\[
    f^*(\tilde{x}) \le \sup f^* \quad \text{and} \quad \ell(A\tilde{x}) \ge \inf \ell,
\]
then, putting these together with the above statements, we find that $\tilde{x}$
and $A\tilde{x}$ are minimizers for the first and second terms, respectively,
\ie,
\[
    f^*(\tilde{x}) = \sup f^* \quad \text{and} \quad \ell(A\tilde{x}) = \inf \ell.
\]
This means that $\tilde{x} \in X^\star$ and $A\tilde{x} \in Y^\star$, or,
equivalently, that $AX^\star \cap Y^\star$ is nonempty. The statement above
follows from the contrapositive: if $AX^\star \cap Y^\star$ is nonempty, then
$p=0$ is not a minimizer for $g$.

\paragraph{Separating hyperplanes.} The prices $p$ that minimize the function
$g$ are intimately related to the geometry of the sets $X^\star$ and $Y^\star$.
(We will see this soon.) For this purpose, we will define $K$ to be the cone of
hyperplanes that separate the sets $AX^\star$ and $Y^\star$, defined:
\[
    K = \{p' \in \reals^m \mid p'^TAx \ge p'^Ty ~\text{for all} ~ x \in X^\star, ~ y \in Y^\star\}.
\]
Note that $0 \in K$, so this set is always nonempty, and $K$ is closed and
convex as it can be written as the intersection of closed halfspaces (which are
also convex sets):
\[
    K = \bigcap_{\substack{x \in X^\star\\y \in Y^\star}} \{p'\in \reals \mid p'^T(Ax - y) \ge 0 \}.
\]

\paragraph{Conditions on prices.}
We will show that, if $p$ satisfies $g(p) \le g(0)$, then $p \in K$. In other
words, any minimizer $p$ of the dual function $g$ must be a separating
hyperplane of the extremizing sets $AX^\star$ and $Y^\star$. The proof is
relatively simple. Since
\[
    g(p) \le g(0) = \sup f^* - \inf \ell,
\]
then, using the definition of $g$, we have that
\[
    f^*(x) - \ell(y) + p^T(y - Ax) \le \sup f^* - \inf \ell,
\]
for every $x \in \conv(S)$ and $y \in \reals^m_+$, so
\[
    p^T(y - Ax) \le \sup f^* - f^*(x) + \ell(y) - \inf \ell.
\]
If we restrict $x$ and $y$ to be in $X^\star$ and $Y^\star$, respectively,
and negate both sides, we then have
\[
    p^T(Ax - y) \ge 0,
\]
or, that $p^TAx \ge p^Ty$ for all $x \in X^\star$ and $y \in Y^\star$, which is
the definition of $p \in K$. (Note that we may replace the inequalities with
strict inequalities and $K$ with $\intr K$, the interior of the cone $K$, to
receive a second useful statement.)

\paragraph{Discussion.} We note that the above definitions serve as a natural
generalization of the condition that the resource utilization is equal to a
target utilization $b^\star$. In our case, we can have many `optimal'
utilizations, given by $Y^\star$, with the additional granularity that any
suboptimal resource utilization vector $y \not \in Y^\star$ has a certain
degree of displeasure, $\ell(y) > \inf \ell$. If the set of optimal utilizations
(for the loss function $\ell$) do not overlap with the zero-cost demand, $X^\star$,
then the original condition states that the prices must be nonzero.

On the other hand, we know that any set of optimal prices $p$ must be a
separating hyperplane for the sets $AX^\star$ and $Y^\star$; \ie, that $p \in
K$. This leads to some interesting observations. If zero utilization is a
possible target, \ie, $0 \in Y^\star$, as is the case for any nondecreasing
loss such as~\eqref{eq:nondec-loss}, then the set $K$ is contained in the dual cone of 
$\cone(AX^\star)$, where $\cone(AX^\star)$ is the set of conic (nonnegative)
combinations of elements in $AX^\star$. For more information, see,
\eg,~\cite[\S2.6]{boyd2004convex}.

\paragraph{Extensions.} There is also a partial converse to the above
conditions on the prices $p$. In particular, for any resource cost $p$ in the
interior of the cone $\intr K$ (satisfying the technical condition that $\ell$
contains no line in the direction of $p$) there is always some scalar $t > 0$
such that $g(tp) < g(0)$; \ie, $0$ cannot be a minimizer for $g$ if the
interior of $K$ is nonempty. While interesting, this point is slightly
technical, so we defer it to appendix~\ref{app:converse}. In general, we can
view this statement as a stronger version of the original claim: if the compact
set $AX^\star$ and closed set $Y^\star$ are disjoint, then there is a strict
separating hyperplane between $AX^\star$ and $Y^\star$, say $p$, so the set
$\intr K$ is nonempty since it contains $p$. This, in turn, would immediately
imply that $g(0)$ cannot be minimal.

\subsection{Properties}\label{sec:properties}
There are a number of properties of the prices $p$ that can be derived from the
dual problem~\eqref{eq:dual-problem}.

\paragraph{Nonnegative prices.}\label{sec:nonnegative} If the objective
function $\ell$ is separable and nondecreasing, as in~\eqref{eq:nondec-loss},
then any price $p_i$ feasible for problem~\eqref{eq:dual-problem} must be
nonnegative, $p_i \ge 0$. (By feasible, we mean that $g(p) < \infty$.) To see
this, note that, by definition~\eqref{eq:nondec-loss}, we have
\[
    \ell^*(p) = \sup_y \left(p^Ty - \ell(y)\right) 
    = \sum_{i=1}^m \sup_{y_i} \left(p_iy_i - \phi_i(y_i)\right),
\]
so we can consider each term individually. If $p_i < 0$ then any $y_i \le 0$ must have
\[
    p_iy_i - \phi_i(y_i) \ge p_iy_i -\phi_i(0) \to \infty,
\]
as $y_i \to -\infty$ since $\phi_i(y_i)$ is nondecreasing in $y_i$. So $g(p) \to
\infty$ and therefore this choice of $p$ cannot be feasible, so we must have
that $p_i \ge 0$.

\paragraph{Superlinear separable losses.}
If the losses $\phi_i$ are superlinear, in that
\begin{equation}\label{eq:superlinearity}
    \frac{\phi_i(z)}{z} \to \infty,
\end{equation}
as $z \to \infty$ and bounded from below, in addition to being nondecreasing,
then $\ell^*(p)$ is finite for $p \ge 0$. This means that the effective domain of
$g$, defined as the set of prices for which $g$ is finite,
\[
    \dom g = \{p \in \reals^m \mid g(p) < \infty\},
\]
is exactly the nonnegative orthant. (This discussion may appear somewhat
theoretical, but we will see later that this turns out to be an important
practical point when updating prices.) While not all losses are superlinear, we
can always make them so by, \eg, adding a small, nonnegative squared term to
$\phi_i$, say
\[
    \tilde\phi_i(z) = \phi_i(z) + \rho(z)_+^2,
\]
where $(z)_+ = \max\{0, z\}$ and $\rho > 0$ is a small positive value, or by
setting $\phi_i(z) = \infty$ for $z \ge 0$ large enough.

\paragraph{Subsidies.} Alternatively, if the function $\ell$ is decreasing
somewhere on the interior of its domain, then there exist points $y^\star$ for
which prices $p_i$ are negative---\ie, sometimes the network is willing to
subsidize usage by paying users to use the network to meet its intended target.
The interpretation is simple: if the base demand of the network is not enough
to meet the target amount, then the network has an incentive to subsidize users
until the marginal cost of the target usage matches the subsidy amounts. 
We note that this may only apply to very specific transaction types in
practice, as it is difficult to issue subsidies in an incentive-compatible
manner that doesn't incentivize the inclusion of `junk' transactions.

\paragraph{Maximum price.} Another observation is that there often exist prices
past which transaction producers would always prefer to not submit a
transaction (or, more generally, will only submit transactions that consume no
resources, if such transactions exist). In fact, we can characterize the set of
all such prices.

To do this, write $S_0\subseteq S$ for the set of transactions bundles that use
no resources, defined
\[
    S_0 = \{x \in S \mid Ax = 0\}.
\]
If $0 \in S$ then $S_0$ is nonempty (as $0 \in S_0$), and we usually expect this
set to be a singleton, $S_0 = \{0\}$. Otherwise, we are saying that there are
transactions that are always costless to include. Now, we will define the set
\[
    P = \bigcap_{x \in S \setminus S_0}\{p \in \reals^m_+\mid p^TAx > q^Tx \}.
\]
This is the set of prices $p \in P$ such that, for every possible transaction
bundle $x \in S$, the price of this transaction bundle, $p^TAx$, paid to the
network, is strictly larger than the total welfare generated by including it,
which is $q^Tx$. (That is, any transaction bundle $x$ that is not costless is
always strictly worse than no transaction, for transaction producers, at these
prices.) The set $P$ is nonempty since $\ones^TAx > 0$ for every $x \in S
\setminus S_0$ (and $S \setminus S_0$ is finite) so, setting $p = t\ones$, we
have that
\[
    p^TAx - q^Tx = t\ones^TAx - q^Tx \to \infty > 0,
\]
as $t \to \infty$, so $t\ones \in P$ for $t$ large enough. The set $P$ is also
a convex set, as it is the intersection of convex sets. Additionally, if $p \in
P$, then any prices $p'$ satisfying $p' \ge p$ must also have $p' \in P$, where
the inequality is taken elementwise. In English: if certain resource prices $p
\in P$ would mean that transactions that consume resources are not included,
then increasing the price of any resources to $p' \ge p$ also implies the same.

\subsection{Solution methods}\label{sec:solution-methods}
As mentioned before, the dual problem~\eqref{eq:dual-problem} is convex. This
means that it can often be easily solved if the function $g$ (or its
subgradients) can be efficiently evaluated. We will see that, assuming users
and miners are approximately solving problem~\eqref{eq:miner-packing}, we can
retrieve approximate (sub)gradients of $g$ and use these to (approximately)
solve the dual problem~\eqref{eq:dual-problem}. In a less-constrained
computational environment, a quasi-Newton method (\eg, L-BFGS) would converge
quickly to the optimal prices and be efficient to implement. However, these
methods aren't amenable to on-chain computation due to their memory and
computational requirements. To solve for the optimal fees on chain, we
therefore propose a modified version of gradient descent which is easy to
compute and does not require additional storage beyond the fees themselves.

\paragraph{Projected gradient descent.} A common algorithm for unconstrained function
minimization problems, such as problem~\eqref{eq:dual-problem}, is \emph{gradient
descent}. In gradient descent, we are given an initial starting point $p^0$
and, if the function $g$ is differentiable, we iteratively update the prices
in the following way:
\[
    p^{k+1} = p^k - \eta \nabla g(p^k).
\]
Here, $\eta > 0$ is some (usually small) positive number referred to as the
`step size' or `learning rate' and $k=0, 1, \dots$ is the iteration number.
This rule has a few important properties. For example, if $\nabla g(p^k) = 0$,
that is, $p^k$ is optimal, then this rule does not update the prices, $p^{k+1}
= p^k$; in other words, any minimizer of $g$ is a fixed point of this update
rule. Additionally, this rule can be shown to converge to the optimal value
under some mild conditions on $g$, \cf~\cite[\S1.2]{bertsekas99}. This update also
has a simple interpretation: if $\nabla g(p^k)$ is not zero, then a small
enough step in the direction of $\nabla g(p^k)$ is guaranteed to evaluate to a
lower value than $p^k$, so an update in this direction decreases the objective
$g$. (This is why the parameter $\eta$ is usually chosen to be small.)

Note that if the effective domain of the function $g$, $\dom g$, is not
$\reals^m$, then it is possible that the $(k+1)$st step ends up outside of the
effective domain, $p^{k+1} \not \in \dom g$, so $g(p^k) = \infty$ which would 
mean that the gradient of $g$ at price $p^{k+1}$ would not exist. To avoid this, we
can instead run \emph{projected} gradient descent, where we project the update
step into the domain of $g$, in order to get $p^{k+1} \in \dom g$, \ie,
\begin{equation}\label{eq:grad-updates}
    p^{k+1} = \proj(p^k - \eta \nabla g(p))
\end{equation}
where $\proj(z)$ is defined
\[
    \proj(z) = \argmin_{p \in \dom g} \|z- p\|_2^2.
\]
In English, $\proj(z)$ is the projection of the price to the nearest point in
the domain of $g$, as measured by the sum-of-squares loss $\|\cdot\|_2^2$.
(This always exists and is unique as the domain of $g$ is always closed and
convex, for any loss function $\ell$ as defined above.) There is relatively
rich literature on the convergence of projected gradient descent, and we refer
the reader to~\cite{shor1985minimization, bertsekas99} for more.

\paragraph{Evaluating the gradient.} In general, since we do not know $q$, we
cannot evaluate the function $g$ at a certain point, say $p^k$. On the other hand,
the gradient of $g$ at $p^k$, when $g$ is differentiable, depends only on the
solution to problem~\eqref{eq:miner-packing} and the maximizer for the dual
function $\ell^*$, at the price $p^k$. (This follows from the gradient equation
in~\eqref{eq:grad}.) So, if we know that transaction producers are solving
their welfare maximization problem~\eqref{eq:miner-packing} to (approximate)
optimality, equation~\eqref{eq:grad} suggests a simple descent algorithm for
solving the dual problem~\eqref{eq:dual-problem}.

From before, let $y^\star$ be a maximizer of $\sup_y \left(y^Tp^k -
\ell(y)\right)$, which is usually easy to compute in practice, and let $x^0$ be
an (approximate) solution to the transaction inclusion
problem~\eqref{eq:miner-packing} (observed, \eg, after the block is built with
resource prices $p^k$). We can approximate the gradient of $g$ at the current
fees $\fees$ using~\eqref{eq:grad}, where we replace the true solution
$x^\star$ with the observed solution $x^0$. Since $x^0$ is Boolean, we can
compute the resource usage $Ax^0$ after observing only the included
transactions. We can then update the fees $\fees^k$ in, say, block $k$, to a
new value $\fees^{k+1}$ by using projected gradient descent with this new
approximation:
\begin{equation}\label{eq:update}
    p^{k+1} = \proj(p^k - \eta (y^\star - Ax^0)).
\end{equation}

\paragraph{Discussion.} 
Whenever $\ell$ is differentiable at $y^\star$, we have that
$\nabla\ell(y^\star) = p$. (To see this, apply the first-order optimality
conditions to the objective in the supremum in the definition of $\ell^\star$.)
We can then think of $y^\star$ as the resource utilization such that the
marginal cost to the network $\nabla\ell$ is equal to the current fees $p$.
Thus, the network aims to set $p$ such that the realized resource utilization
is equal to $y^\star$. We can see that~\eqref{eq:grad-updates} will increase
the network fee for a resource being overutilized and decrease the network fee
for a resource being underutilized. This pricing mechanism updates fees to
disincentivize future users and miners from including transactions that consume
currently-overutilized resources in future blocks. Additionally, we
note that algorithms of this form are not the only algorithms which are reasonable.
For example, any algorithm that has a fixed point $p$ satisfying $\nabla g(p) = 0$
and converges to this point under suitable conditions is also similarly reasonable. 
One well-known example is an update rule of the form of~\eqref{eq:buterinUpdate}:
\[
    p^{k+1} = p^k \odot \exp(-\eta \nabla g(p^k)),
\]
when the prices must be nonnegative, \ie, when $\dom g \subseteq \reals_+^m$.
We note that one important part of reasonable rules is that they only depend on
(an approximation of) the gradient of the function $g$, since the value of $g$
may not even be known in practice. Additionally, in some cases, the function
$g$ is nondifferentiable at prices $p$. In this case, the subgradient still
often exists and convergence of the update rule~\eqref{eq:update} can be
guaranteed under slightly stronger conditions. (The modification is needed as
not all subgradients are descent directions.)

\paragraph{Simple examples.}
We can derive specific update rules by choosing specific loss functions.
For example, consider the loss function
\[
    \ell(y) = \begin{cases}
        0 & y = b^\star \\
        +\infty & \text{otherwise},
    \end{cases}
\] 
which captures infinite unhappiness of the network designer for any deviation
from the target resource usage $b^\star$. The corresponding conjugate function
is
\[
    \ell^*(p) = \sup_y(y^Tp - \ell(y)) = (b^\star)^Tp,
\]
with maximizer $y^\star = b^\star$. (Note that this maximizer
does not change for any price $p$). Since $\dom g = \reals^m$, then the
update rule is
\[
    \fees^{k+1} = \fees^k - \eta (b^\star - Ax^0).
\]
If the utilization $Ax^0$ lies far below $b^\star$, the fees $p^k$ might become
negative, \ie, the network would want to subsidize certain resource usage to
meet the requirement that it must be equal to $b^\star$.

\paragraph{Nondecreasing separable loss.} A more reasonable family of loss
functions would be the nondecreasing, separable losses:
\[
    \ell(y) = \sum_{i=1}^m \phi_i(y_i).
\]
From~\S\ref{sec:nonnegative} we know that the domain of $g$ is precisely the
nonnegative orthant when the functions $\phi_i$ are superlinear (\ie,
satisfy~\eqref{eq:superlinearity}) and bounded from below, so we have that
\[
    \proj(z) = (z)_+
\]
where $(w)_+ = \max\{0, w\}$ for scalar $w$ and is applied elementwise for vectors.
Additionally, using the definition of the separable loss, we can write
\[
    \ell^*(p) = \sum_{i=1}^m \sup_{y_i} \left(y_ip_i - \phi_i(y_i)\right).
\]
Letting $y^\star_i$ be the maximizers for the individual problems at the current
price $p^k$, we have
\[
    p^{k+1} = (p^k - \eta (y^\star - Ax^0))_+.
\]
For example, if $\phi_i$ is an indicator function with $\phi_i(y_i) = 0$ if $y_i
\le b_i^\star$ and $\infty$ otherwise, as in the
loss~\eqref{eq:inequality-loss}, then an optimal point is always $y_i^\star =
b_i^\star$, when $p_i \ge 0$. Since no updates will ever set $p_i < 0$, we
therefore have,
\[
    p^{k+1} = (p^k - \eta (b^\star - Ax^0))_+.
\]
which is precisely the update given in~\eqref{eq:basic-update}.
In addition, a simple transform applied to the price variable gives an update 
rule similar to~\eqref{eq:buterinUpdate} (see appendix~\ref{app:exp-update}).

\paragraph{Notes.}
While we used projected gradient descent rules for the examples above, we note
that this class of update rules is not the only option. Other update rules
naturally fall out of other optimization algorithms applied
to~\eqref{eq:dual-problem}. For example, if we only want to update some subset
of the prices at each iteration, we can use block coordinate descent. We can
also add adaptive step size rules or momentum terms to our gradient descent
formulation. These additions would yield more computationally expensive
algorithms, but they may result in faster convergence to optimal prices when
the distribution of processed transactions shifts. This is a potentially
interesting area for future research.

\section{The cost of uniform prices}
In this section, we show that pricing resources using the method outlined above
can increase network efficiency and 
make the network more robust to DoS attacks or distribution shifts.
We construct a toy experiment to illustrate these differences between uniform 
and multidimensional resource pricing, and leave more extensive numerical studies to future work.

\paragraph{The setup.}
We consider a blockchain system with two resources (\eg, compute and storage) with 
resource utilizations $y_1$ and $y_2$.
Resource 1 is much cheaper for the network to use than resource 2, so
\[
    b^\star = \bmat{10 \\ 1} \qquad
    \text{and} \qquad \bmat{y_1 \\ y_2} \leq b = \bmat{50 \\ 5}.
\]
Furthermore, we assume that there is a joint capacity constraint on these resources
\[
    y_1 + 10y_2 \leq 50,
\]
which captures the resource tradeoff.
Each transaction $a_j$ is therefore a vector in $\reals^3_+$ with 
\[
    a_j = \bmat{a_{1j} \\ a_{2j} \\ a_{1j} + 10a_{2j}}.
\]
As in \S\ref{sec:solution-methods}, we consider the simple loss function
\[
    \ell(y) = \begin{cases}
        0 & y = b^\star \\
        +\infty & \text{otherwise},
    \end{cases}
\]
which has update rule
\[
    p^{k+1} = p^k - \eta (b^\star - Ax^0).
\]
In the scenarios below, we compare our multidimensional fee market approach to
a baseline, where both resources are combined into one equal to $a_{1j} + 10a_{2j}$ with 
$b^\star = 20\% \times \max(b_1, b_2) = 10$. 
We demonstrate that pricing these resources separately leads to better network performance.
All code is available at
\begin{center}
    \texttt{https://github.com/bcc-research/resource-pricing}.
\end{center}
We run simulations in the Julia programming language~\cite{julia}.
The transaction producers' optimization problem~\eqref{eq:miner-packing} is modeled with
JuMP~\cite{jump} and solved with COIN-OR's simplex-based 
linear programming solver, Clp~\cite{clp-solver}.
The solution is usually integral, but when it is not, we fall back to the HiGHS
mixed-integer linear program solver~\cite{huangfu2018highsg}.

\paragraph{Scenario 1: steady state behavior.}
We consider a sequence of $250$ blocks. At each block, there are $15$ submitted transactions,
with resource usage randomly drawn as $a_{1j} \sim \mathcal{U}(0.5,1)$ and 
$a_{2j} \sim \mathcal{U}(0.05,0.1)$.
(For example, these may be moderate compute and low storage transactions.)
Transaction utility is drawn as $q_j \sim \mathcal{U}(0, 5)$.
We initialize the price vector as $p = 0$ and examine the steady state behavior,
where the price updates and transaction producer behavior are defined as in the previous section.
We use a learning rate $\eta = 1\times 10^{-2}$ throughout.

\begin{figure}
    \centering
    \includegraphics[width=0.48\linewidth]{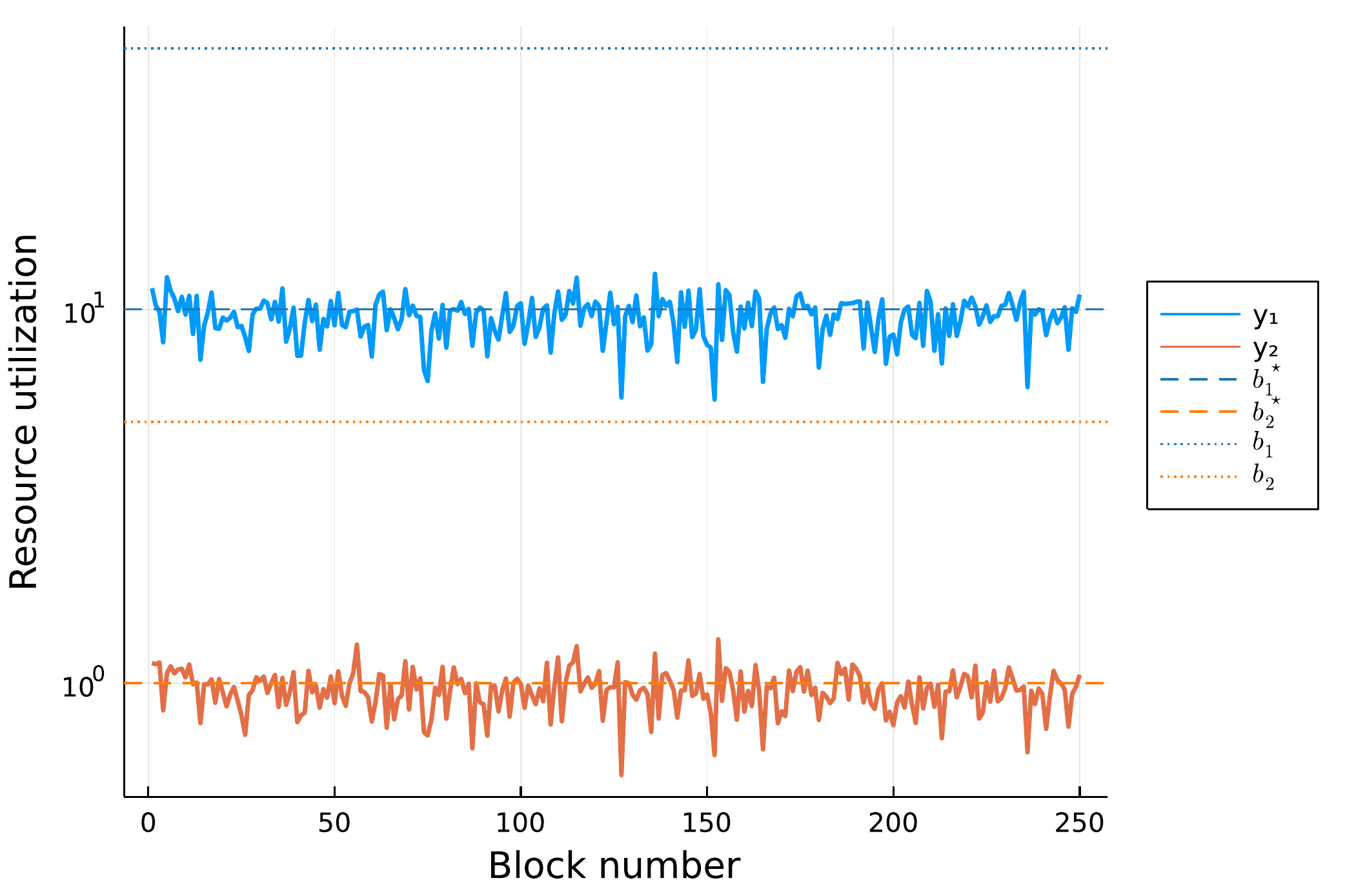}
    \includegraphics[width=0.48\linewidth]{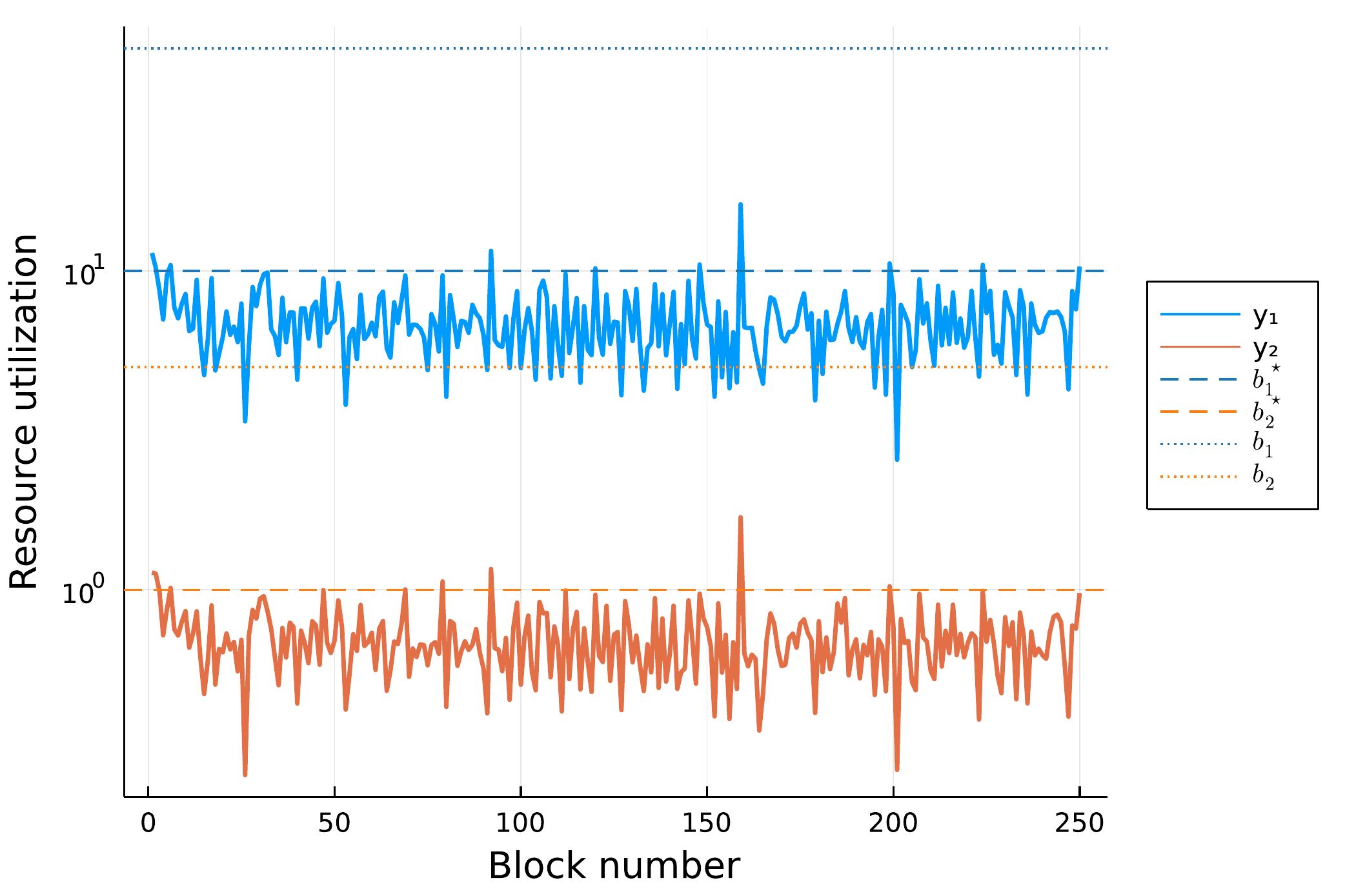}
    \caption{
        Resource utilization for multidimensional pricing (left) clusters
        closer to target values than for uniform pricing (right), which includes
        limited information about individual targets or caps.
    }
    \label{fig:scenario1-rx}
\end{figure}

\begin{figure}
    \centering
    \includegraphics[width=0.48\linewidth]{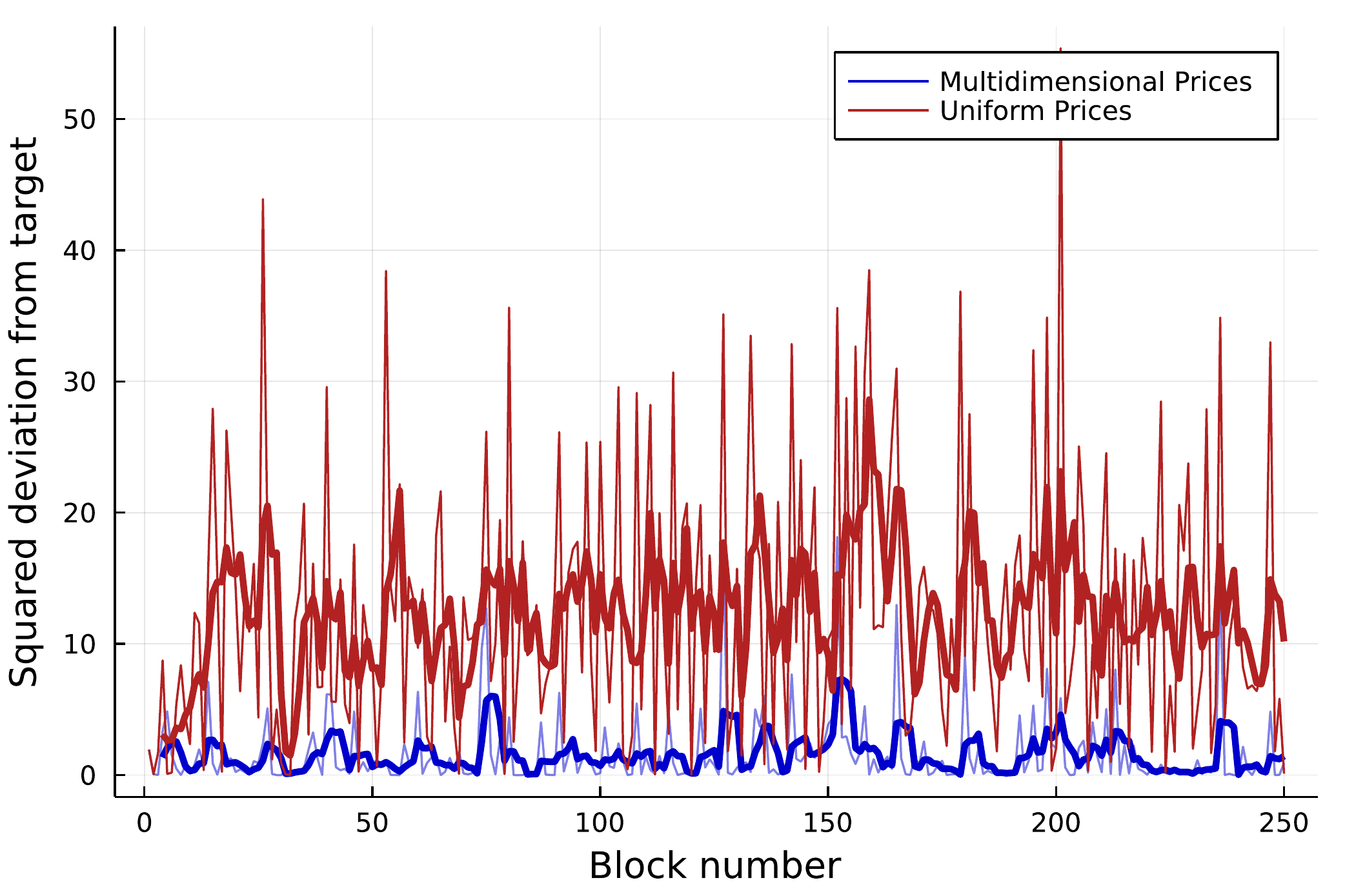}
    \includegraphics[width=0.48\linewidth]{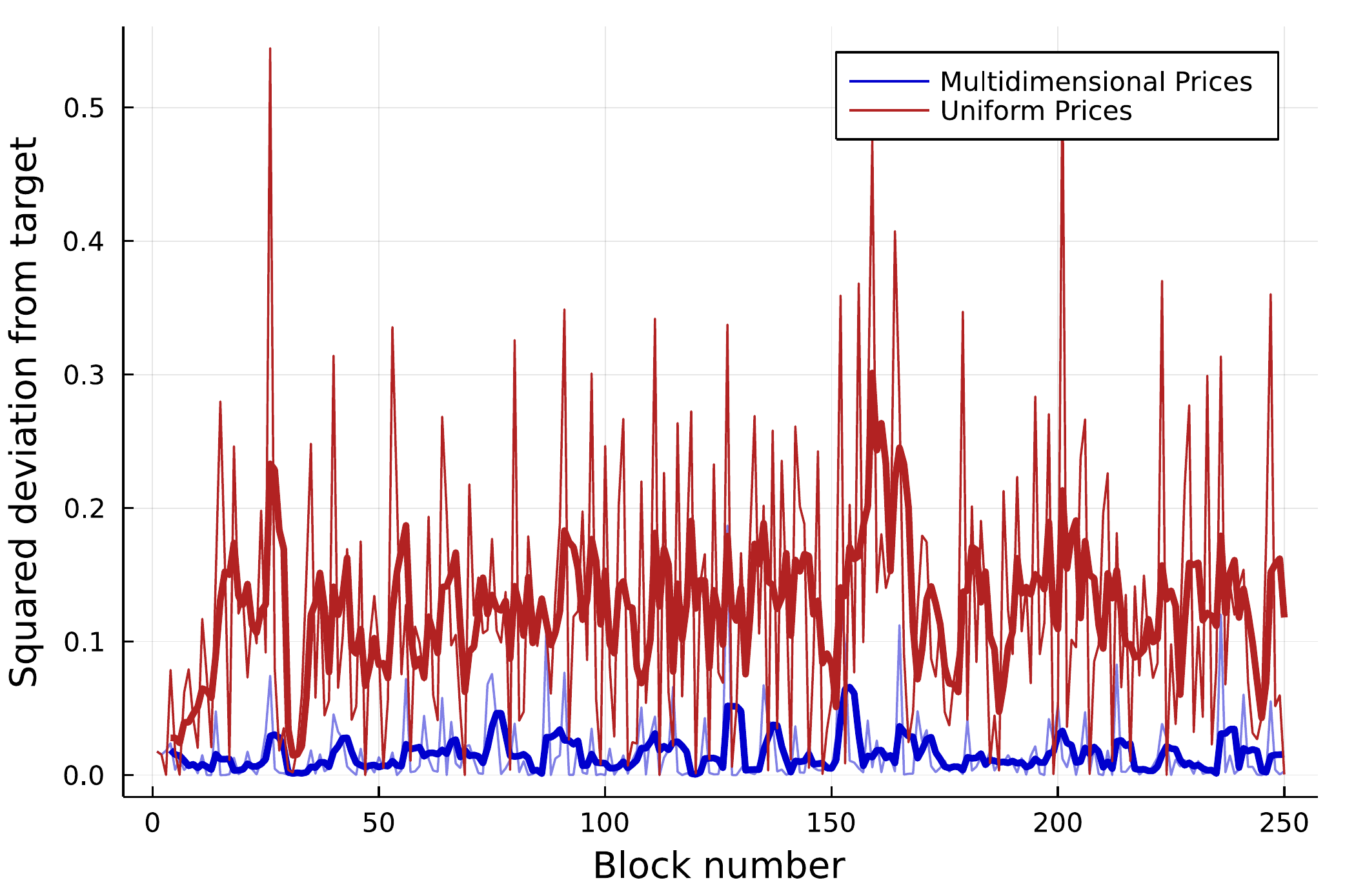}
    \caption{
        Squared deviation from the target values given by $b^\star$ in
        resource 1 utilization (left) and resource 2 utilization (right)
        for both uniform and multidimensional pricing.
        The thicker line is the four-sample moving average.
    }
    \label{fig:scenario1-dev}
\end{figure}

\begin{figure}[!ht]
    \centering
    \includegraphics[width=0.48\linewidth]{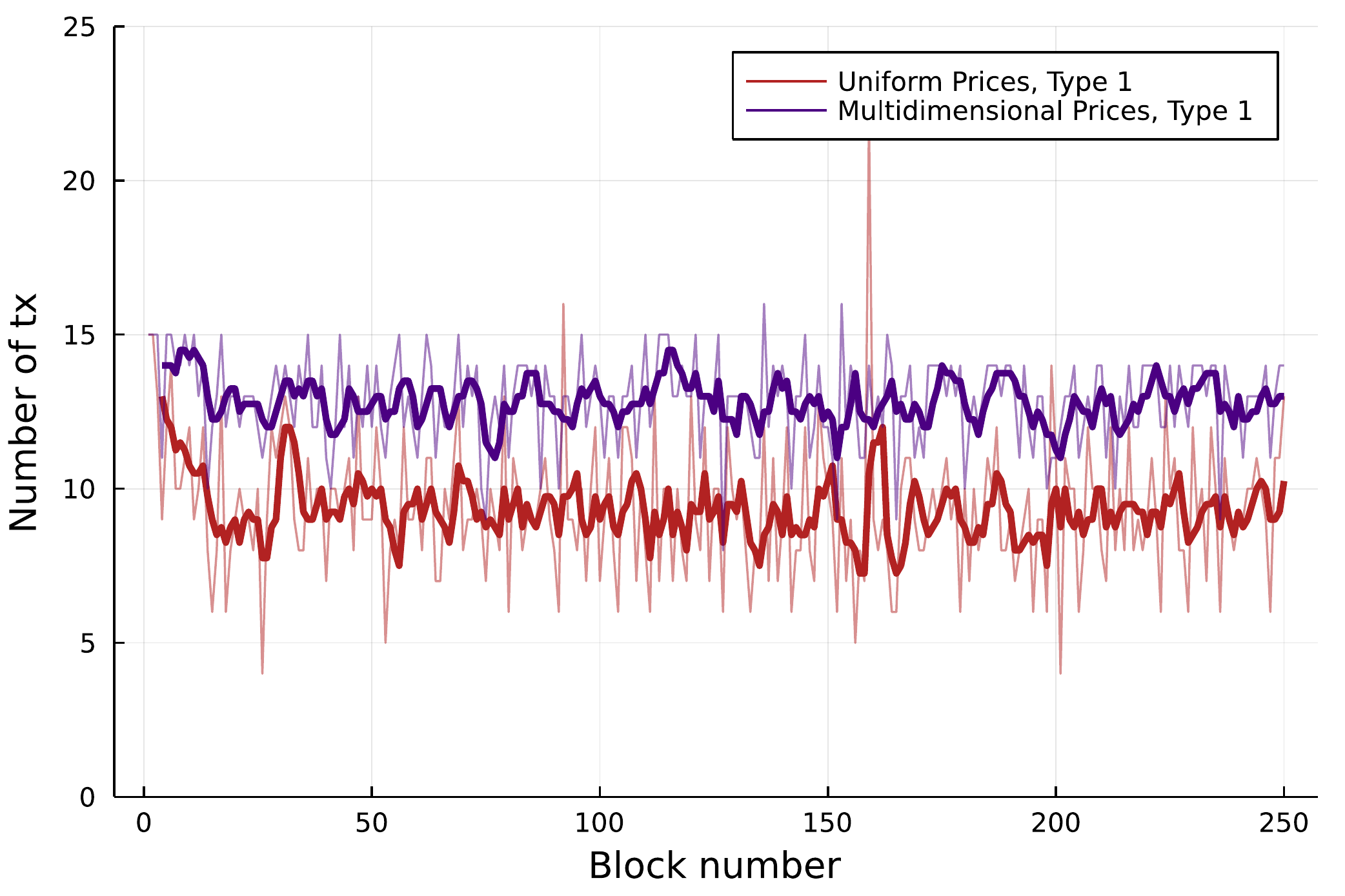}
    \includegraphics[width=0.48\linewidth]{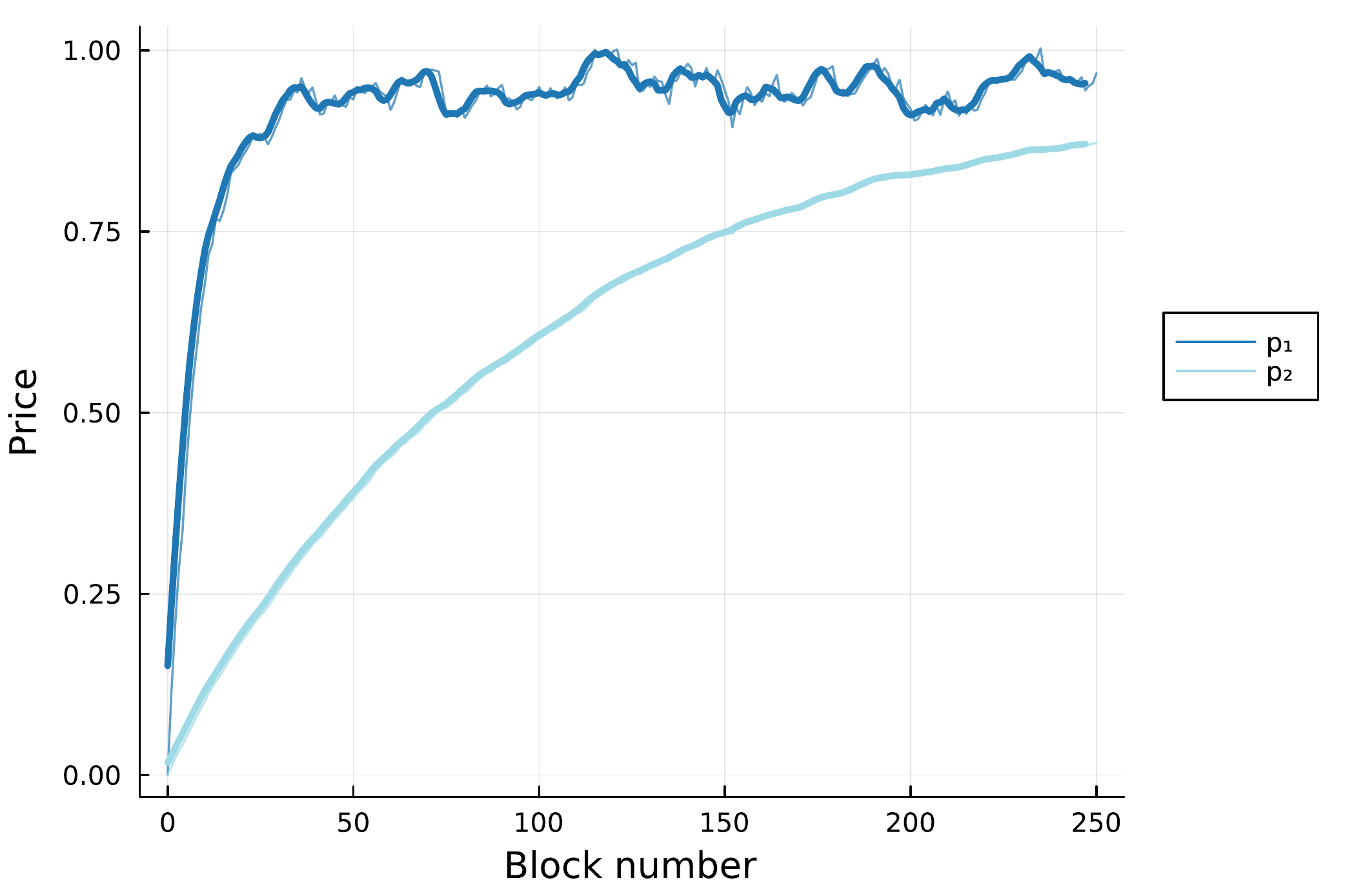}
    \caption{
        Multidimensional pricing allows us to include more transactions per block (left)
        by optimally adjusting prices (right). The thicker line is the four-sample
        moving average of the data.
    }
    \label{fig:scenario1-ntx}
\end{figure}

The resource utilization, shown in figure~\ref{fig:scenario1-rx}, suggests that our multidimensional
scheme more closely tracks the target utilities $b^\star$ than a single-dimensional
fee market. 
Figure~\ref{fig:scenario1-dev} shows the squared deviation from the target 
resource utilizations.
Furthermore, the number of transactions included per block is consistently
higher, illustrated in figure~\ref{fig:scenario1-ntx}
(purple line).

\paragraph{Scenario 2: transactions distribution shift.}
Often, the distribution of transaction types submitted to a blockchain network 
differs for a short period of time (\eg, during NFT mints).
There may be a change in both the number of transactions and the distribution
of resources required. We repeat the above simulation but add $150$
transactions in block $10$; each transaction has a resource vector $a_j = (0.01, 0.5)$.
(For example, these transactions may have low computation but high storage requirements.)
We draw the utility $q_j \sim \mathcal{U}(10, 20)$ and begin the network at the
steady-state prices from scenario 1.

\begin{figure}
    \centering
    \includegraphics[width=0.48\linewidth]{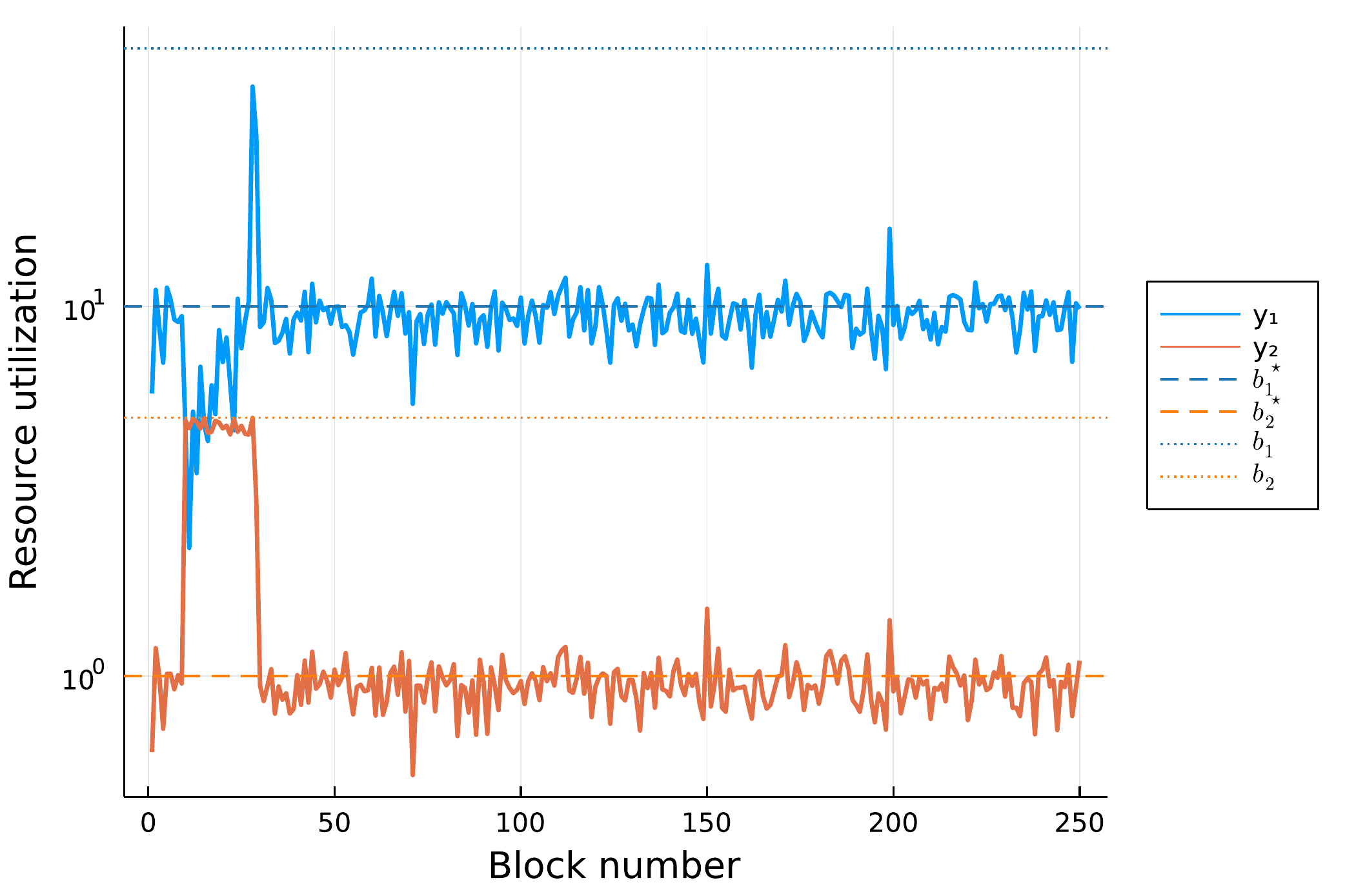}
    \includegraphics[width=0.48\linewidth]{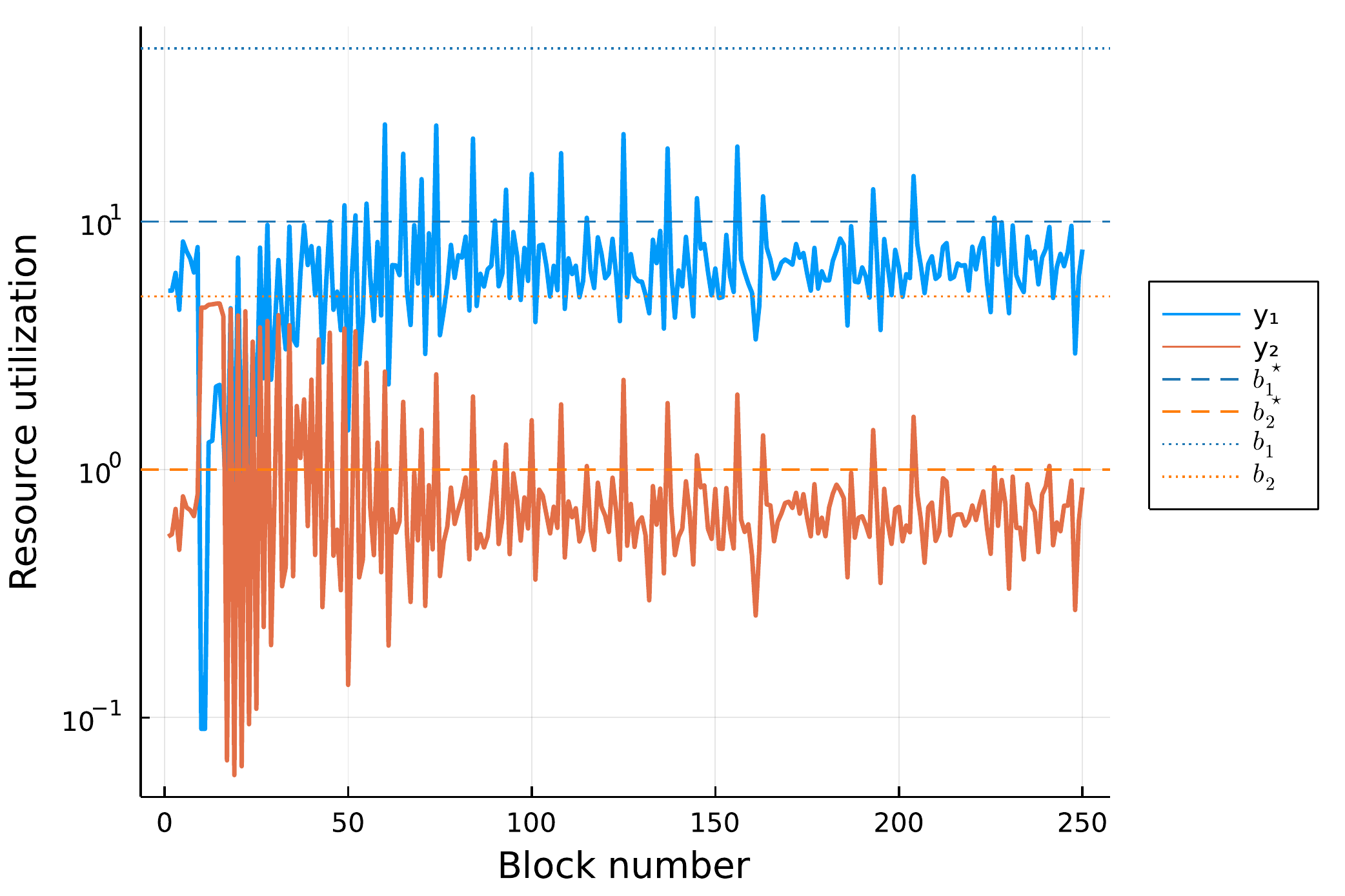}
    \caption{
        Resource utilization for multidimensional pricing (left) clusters
        closer to target values than for uniform pricing (right) after a burst
        to the resource limit to handle transactions that make heavy use of 
        resource 2.
    }
    \label{fig:scenario2-rx}
\end{figure}

\begin{figure}[!ht]
    \centering
    \includegraphics[width=0.48\linewidth]{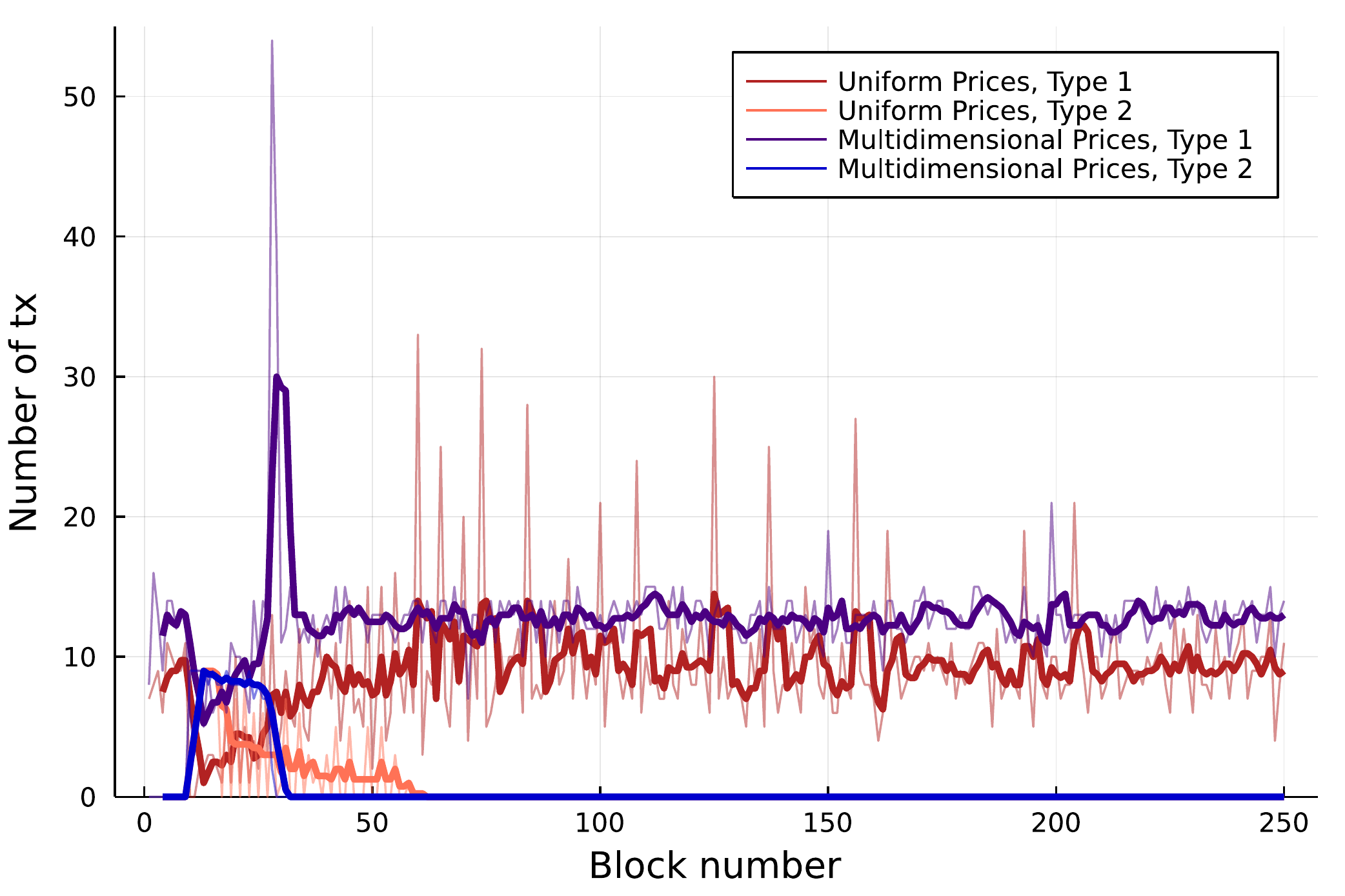}
    \includegraphics[width=0.48\linewidth]{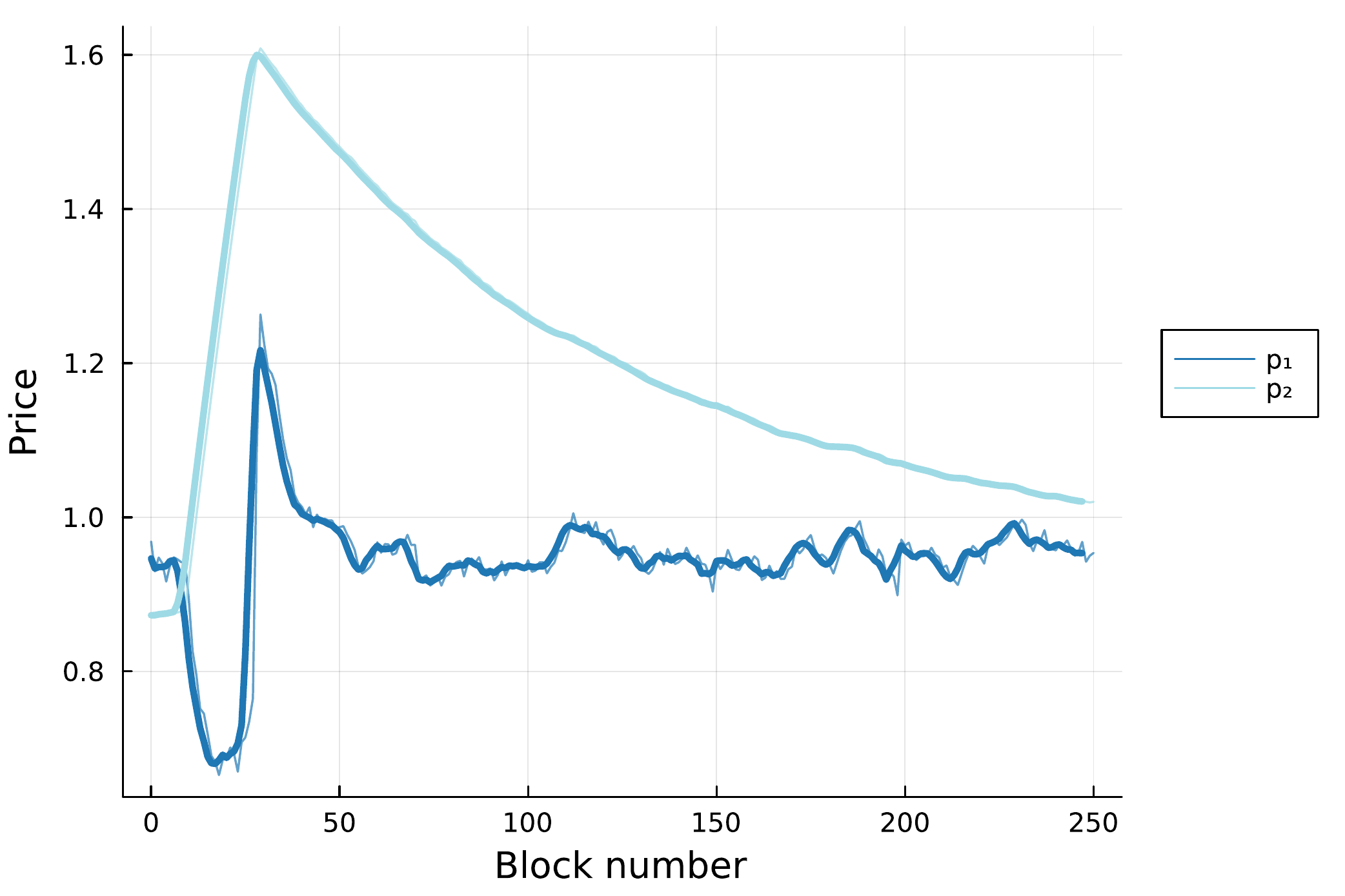}
    \caption{
        Multidimensional pricing allows us to include more transactions per
        block by optimally adjusting prices (right). The thicker line is the
        four-sample moving average of the data.
    }
    \label{fig:scenario2-ntx}
\end{figure}

In figure~\ref{fig:scenario2-rx}, we see that a multidimensional fee market gracefully
handles the distribution shift.
The network fully utilizes resource 2 for a short period of time before 
returning to steady state.
Uniform pricing, on the other hand, does not do a good job of adjusting its
resource usage and oscillates around the target.
Figure~\ref{fig:scenario2-ntx} show that, as a result, multidimensional pricing
is able to include more transactions, both during the distribution shift
and after the network returns to steady state. 
We see that the prices smoothly adjust accordingly.

\section{Extensions}\label{sec:Extensions}
\paragraph{Parallel transaction execution model.}
\label{app:parallel}
Consider the scenario where the nodes have $L$ parallel execution environments
(\eg, threads), 
each of which has its own set of $m$ identical resources.
In addition, there are $r$ resources shared between the environments.
We denote transactions run on thread $k$ by $x_k \in \{0,1\}^n$.
The resource allocation problem becomes
\[
    \begin{aligned}
        &\text{maximize}     && \sum_{k=1}^L q^Tx_k - \ell(y_1, \dots, y_L, y^\mathrm{shared}) \\
        &\text{subject to}   && y_k = Ax_k,          \qquad k = 1, \dots, L      \\
                            &&& z = \sum_{k=1}^L x_k \\
                            &&& y^\mathrm{shared} = Bz \\
                            &&& z \in \conv(S^\mathrm{shared}) \\
                            &&& x_k \in \conv(S), \qquad k = 1, \dots, L.
    \end{aligned}
\]
As before, the Boolean vector sets $S$ and $S^\mathrm{shared}$ encode
constraints such as resource limits for each parallel environment and the
shared environment respectively. In addition, we'd expect to have $z \leq
\ones$ if each transaction is only allocated to a single environment, which can
be encoded by $S^\mathrm{shared}$. By stacking the variables into one vector,
this problem can be seen as a special case of~\eqref{eq:allocation} and can be
solved with the same method presented in this work. (The interpretation here
is that we are declaring a number of `combined resources', each corresponding
to the parallel execution environments along with their shared resources.)

\paragraph{Different price update speeds.} Some resources may be able to
sustain burst capacities for much shorter periods of times than other
resources. In practice, we may wish to increase the prices of these resources
faster.
For example, a storage opcode that generates a lot of memory allocations will 
quickly cause garbage collection overhead, which could slow down the network.
As a result, we likely want to increase its price faster than the prices of basic 
arithmetic, even under the same relative utilization.
To do this, we can update~\eqref{eq:grad-updates} to include a learning rate for 
each resource.
We collect these in a diagonal matrix $D = \diag{\eta_1, \dots, \eta_m}$:
\[
    p^{k+1} = \proj(p^k - D\nabla g(p))
\]
These learning rates can be chosen by system designers using simulations and historical data.

\paragraph{Contract throughput.}
Alternatively, we can define utilization on a per-contract basis 
instead of a per-resource basis (per-contract fees were recently proposed by the
developers of Solana~\cite{aeyakovenko_consider_2021}).
We define the utilization of a smart contract $j$ as $z_j = (w^Ta_j)x_j$,
where $w$ is some weight vector and $x_j \in \integers_+$ is the number of times 
contract $j$ is called.
In matrix form, $z = A^Tw \odot x$, where $\odot$ is the Hadamard (elementwise) product.
For each contract, the utilization $z_j$ is $0$ when $x_j = 0$, 
which can be interpreted as not calling contract $j$ in a block.
When $x_j > 0$, the utilization is $\left(\sum_i w_i a_{ij}\right)x_j$.
When we use per-contract utilizations, the loss function can capture a notion 
of fairness in resource allocation to contracts.
For example, we may want to prioritize cheaper-to-execute contracts over more 
expensive ones by using, \eg, proportional fairness as in~\cite{kelly1997charging},
though there are many other notions that may be useful.
With this setup, the resource allocation problem is
\[
    \begin{array}{ll}
        \text{maximize}     & q^Tx - \ell(z) \\
        \text{subject to}   & z = A^Tw \odot x \\
                            & x \in \conv(S).
    \end{array}
\]
Again, we can introduce the dual variable $\fees \in \reals^n$ for the equality constraint,
and, with a similar method to the one introduced in this paper, iteratively update this 
variable to find the optimal fees to charge for each smart contract call.

\section{Conclusion}
We constructed a framework for multidimensional resource pricing in blockchains.
Using this framework, we modeled the network designer's goal of maximizing
transaction producer utility, minus the loss incurred by the network, as a
an optimization problem.
We used tools from convex optimization---and, in particular, duality theory---to
decompose this problem into two simpler problems: one solved on chain by
the network, and another solved off chain by the transaction producers.
The prices that unify the competing objectives of minimizing network loss 
and maximizing transaction producer utility are precisely the dual variables
in the optimization problem.
Setting these prices correctly (\ie, to minimize the dual function) results in
a solution to the original problem.
We then demonstrated efficient methods for updating prices that are amenable to on-chain computation.
Finally, we numerically illustrate, via a simple example, the proposed pricing mechanism.
We find that it allows the network to equilibrate to its resource utilization
target more quickly than the uniform price case, while offering greater
throughput without increasing node hardware requirements.

To the best of the authors' knowledge, this is the first work to systematically
study optimal pricing of resources in blockchains in the many-asset setting.
Future work and improvements to this model include a detailed
game-theoretic analysis, extending that of~\cite{ferreira2021dynamic}, along
with a more concrete analysis of the dynamical behavior of fees set in
this manner.
Finally, a more thorough numerical evaluation of these methods under realistic
conditions (such as testnets) will be necessary to see if these methods are
feasible in production.

\section*{Acknowledgements}
We would like to thank John Adler, Vitalik Buterin, Dev Ojha, Kshitij Kulkarni, 
Matheus Ferreira, Barnab\'{e} Monnot, and Dinesh Pinto 
for helpful conversations, insights, and edits.
We're especially appreciative to John Adler for bearing with us through many
drafts of this work and consistently providing valuable feedback.

\printbibliography

\appendix
\section{A (very short) primer on convexity} \label{app:convex}
This appendix serves as a short introduction to the basic notions of convexity
used in this paper for readers familiar with basic real analysis and linear
algebra. For (much) more, we recommend~\cite{boyd2004convex}.

\subsection{Basic definitions}
\paragraph{Convexity.} We say a set $S\subseteq \reals^m$ is convex if, for any two points $x, y \in S$
and any $0 \le \gamma \le 1$, we have
\[
\gamma x + (1-\gamma) y \in S.
\]
In other words, a set $S$ is convex if it contains all line segments between
any two points in $S$.
A classic example of a closed convex set is a closed halfspace
\[
    H = \{x \in \reals^m \mid a^Tx \le b\},
\]
where $a \in \reals^m$ and $b \in \reals$. (An otherwise silly but useful
example is the empty set, which vacuously meets the requirements.)

We say a function over the extended reals $f: S \to \reals \cup
\{\infty\}$ is convex if $S$ is convex and, for any $x, y \in S$ and $0 \le
\gamma \le 1$,
\[
    f(\gamma x + (1-\gamma)y) \le \gamma f(x) + (1-\gamma) f(y).
\]
Equivalently: a function is convex if any chord of the function (a line
segment between two points on its graph) lies above (or, strictly speaking, not below) the
function itself. We say $f$ is concave if $-f$ is convex. Some basic functions
that are convex are linear functions $f(x) =q^Tx$ for some $q \in \reals^n$,
norms $f(x) = \|x\|$, and indicator functions of convex sets:
\[
    f(x) = \begin{cases}
        0 & x \in S\\
        \infty & \text{otherwise},
    \end{cases}
\]
where $S \subseteq \reals^m$ is a convex set.

\paragraph{Domain.} Usually it is simpler to work with functions $f$ defined
over all of $\reals^m$ as opposed to just a subset. We may extend the functions
$f: S \to \reals \cup \{\infty\}$ to a function $\tilde f: \reals^m \to \reals
\cup \{\infty\}$ by setting $\tilde f(x) = f(x)$ for $x \in S$ and $\tilde f(x)
= \infty$ if $x \not \in S$. We write the \emph{effective domain} of a function
$f$ defined over $\reals^m$ as
\[
    \dom f = \{x \in \reals^m \mid f(x) < \infty\}.
\]
Throughout the paper and the remainder of this appendix, we assume
that all functions are extended in this way.

\paragraph{Characterizations of convexity.} There are many
equivalent characterizations of convexity for functions. A particularly useful
one, if the function $f$ is differentiable over its domain $\dom f$, is, for any two points $x,
y \in \dom f$, we have:
\begin{equation}\label{eq:conv-ineq}
    f(y) \ge f(x) + \nabla f(x)^T(y-x).
\end{equation}
In other words, any tangent plane to $f$ at $x$ is a global underestimator of
the function. A common characterization for twice-differentiable functions
$f$ is that the hessian (the matrix of all second derivatives) of $f$ at every
point is positive semidefinite, but this is often particularly restrictive.
The gradient-based definition of convexity immediately implies that
if we find a point $x$ with $\nabla f(x) = 0$, then
\[
    f(y) \ge f(x),
\]
for any $y \in \dom f$; \ie, $x$ is a global minimizer of $f$. (The converse
is similarly easy to show.)

\paragraph{Consequences of convexity.} There are a number of important
consequences of convexity. The simplest is: given two closed convex sets
$S, T \subseteq \reals^n$ with $S \cap T = \emptyset$ then there exists
a vector $p \in \reals^n$ and $p \ne 0$ that separates these sets, \ie,
\[
    p^T(x - y) \ge 0, ~ \text{for all} ~ x \in S, ~ y \in T.
\]
If one of the sets, say $S$, is also compact, then it is possible to make the
stronger claim that there exists $p' \in \reals^m$ and $\eps > 0$ such that
\[
    p^T(x - y) \ge \eps, ~ \text{for all} ~ x \in S, ~ y \in T.
\]
The arguments for each of these are relatively simple and can be found
in~\cite[\S2.5]{boyd2004convex}. Nearly all major results in convex
optimization theory are a consequence of these two facts.

\paragraph{Convexity-preserving operations on sets.} There are a number of
operations which preserve convexity of sets. For example, any (potentially
uncountable) intersection of convex sets is convex. The (finite) sum of convex
sets $S, T \subseteq \reals^m$, defined
\[
    S+T = \{x+y \mid x \in S, ~ y \in T\},
\]
is convex, while negation of a set $-S = \{-x \mid x \in S\}$ is also convex.
Any linear function of a convex set, say $A \in \reals^{n\times m}$ with $AS =
\{Ax \mid x \in S\}$, is convex. In the special case that $S$ is also compact,
then $AS$ is compact. (It is, on the other hand, not true in general that if
$S$ is closed then $AS$ is closed.) All of these conditions can be easily
verified from the definitions above. Additionally, there are other operations
that preserve convexity, such as the perspective transform~\cite[\S
2.3.3]{boyd2004convex}, but we do not need those here.

\paragraph{Convexity-preserving operations on functions.} Similar to the
case with convex sets, 
there are a number of convexity-preserving operations on convex functions.
The sum of convex functions is convex, while any nonnegative scaling $\gamma \ge 0$
of a convex function $f$ is convex, \ie, $\gamma f$ is convex. Affine precomposition
of convex functions is convex, \ie, if $f(x)$ is convex over $x \in \reals^m$ then $f(Ay + b)$
is convex over $y \in \reals^n$, for any $A \in \reals^{m \times n}$ and $b \in \reals^m$.
Convexity is preserved over suprema:
\[
    f(x) = \sup_{y \in Y} f_y(x)
\]
where $f_y$ is a family of convex functions indexed by some (potentially
uncountable) set $Y$. The set $Y$ has no assumptions on it, other than being
nonempty. A classic example is if
\[
    f(x) = \sup_{y \in Y} \,\left(y^Tx - g(y)\right),
\]
where $g: Y \to \reals \cup \{\infty\}$ is any (potentially nonconvex) function
and $Y\subseteq \reals^m$ is any nonempty set, then the function $f$ is convex
as it is the supremum of a family of affine (and therefore convex) functions of
$x$. (The function $f$ is known as the \emph{Fenchel conjugate}, often just the
conjugate, of $g$ and is denoted $g^*(x)$.) 
As before, all of these statements are easy to show given
the definitions provided above.

\paragraph{Lines.} In some cases it is important that convex functions
have at least a certain amount of growth. To this end, we say the convex
function $f$ \emph{contains a line} $p \in \reals^m$ with $p \ne 0$ if, for
some $x_0 \in \reals^m$ and $q \in \reals$, we have
\[
    f(x_0 + tp) \le f(x_0) + tq,
\]
for every $t \ge 0$. We say that $f$ \emph{contains a line in the direction of}
$r \in \reals^m$ if it contains a line $p$ with $r^Tp > 0$. It is not hard to show
that functions which satisfy
\[
    \frac{f(tp)}{t} \to \infty
\]
as $t \to \infty$ for every $p \ne 0$ contain no lines. Examples of such functions
are set indicator functions for compact sets and square norms, $f(x) = \|x\|^2$.

\paragraph{Sublevel sets and semicontinuity.} The $\alpha$-sublevel set of a
function $f$ is defined as
\[
    S_\alpha = \{x \in \reals^m \mid f(x) \le \alpha\}.
\]
If the function $f$ is convex, the set $S_\alpha$ is also convex for every
$\alpha$. If the sublevel sets $S_\alpha$ are closed for every $\alpha$, then
we say the function $f$ is \emph{lower semicontinuous} (some authors say the
function $f$ is closed, instead). Any continuous function is lower
semicontinuous, and any lower semicontinuous function $f$ minimized over a set
$S$ such that $\dom f \cap S$ is compact achieves its minimum within that set
\[
    \inf_{x \in S} f(x) = f(x^\star),
\]
for some $x^\star \in S$.

\subsection{Convex hulls}\label{app:convex-hull}
Sometimes the sets that we deal with are not convex, but, in some special
cases, we may treat them as if they were.

\paragraph{Convex hull.} To this end, define the \emph{convex hull} of a set $S
\subseteq \reals^m$, written $\conv(S)$, as the set containing all convex
combinations of points in $S$. Here, a \emph{convex combination} of points in
$S$, say $x_i \in S$ for $i=1, \dots, n$, is any point $y \in \reals^m$
(potentially not in $S$) which can be written as 
\[
    y = \gamma_1x_1 + \dots + \gamma_nx_n,
\]
where $\gamma_i \ge 0$ and $\sum_{i=1}^n \gamma_i = 1$. The set $\conv(S)$ is
evidently convex (we encourage the reader to show this) and is, in fact, the
smallest convex set which contains $S$. If the set $S$ is compact, then its
convex hull is also compact.

\paragraph{Linear functions.} One special case where we may replace a set $S$
with its convex hull is the maximization of linear functions. That is,
for $S \subseteq \reals^m$ and $q \in \reals^m$, we have
\[
    \sup_{x \in S} \,q^Tx = \sup_{x \in \conv(S)} q^Tx
\]
To see this, we will show that any point $x \in \conv(S)$ with value $q^Tx$ has
a point $y \in S$ with $q^Ty \ge q^Tx$. Since $x \in \conv(S)$, then it can be
written as
\[
    x = \gamma_1 x_1 + \dots + \gamma_n x_n,
\]
where $x_i \in S$ and $\gamma_i \ge 0$ with $\sum_{i=1}^n \gamma_i = 1$. So we
have that
\[
    q^Tx = \gamma_1 q^Tx_1 + \dots + \gamma_n q^Tx_n \le (\gamma_1 + \dots + \gamma_n)q^Tx_j =q^Tx_j,
\]
where $j$ is the index for which $q^Tx_j$ is the largest, \ie, $q^Tx_j \ge
q^Tx_i$ for $i=1, \dots, n$. Since $x_j \in S$, then we are done. (This proof
can be extended from linear functions to quasiconvex functions nearly
immediately, but we only use this special case here.)
The other direction immediately follows from the fact that $S \subseteq \conv(S)$.
We note that while the maximum values of these two problems are the same, 
the set of maximizers may not be.

\paragraph{Derivatives.} If the set $S$ is finite, then the function
\[
    g(p) = \max_{x \in S} p^Tx
\]
is differentiable at points $p$ when the optimal point $x^\star \in S$ is
unique, and its derivative at $p$ is $\nabla g(p) = x^\star$. This is easy to see, as $x^\star$ is
unique only when
\[
    p^Tx^\star > p^Tx, \quad x \in S \setminus \{x^\star\},
\]
so $g(p') = p'^Tx^\star$ for $p'$ in a neighborhood of $p$. Differentiating
both sides gives the result.

\subsection{Cones}
A special case of convex sets are the convex cones. A set $K \subseteq \reals^m$ is a \emph{cone}
if $x \in K$ means that $\gamma x \in K$ for any $\gamma \ge 0$. A cone $K$ is a \emph{convex cone}
if $K$ is also convex, which means that it is closed under conic (nonnegative) combinations:
\ie, for any $\gamma_i \ge 0$ and $x_i \in K$ for $i=1, \dots, n$, we have
\[
\gamma_1 x_1 + \dots + \gamma_n x_n \in K.
\]
Convex cones are generally useful as they imply a partial ordering with respect
to their elements, given by $x \gek y$ if $x - y \in K$, and this ordering is
closed under nonnegative multiplications and additions of inequalities.
See~\cite[\S 2.4 \& 2.6]{boyd2004convex} for more information on cones and conic duality.

\subsection{Convex optimization problems}
In many cases, we are interested in optimization problems over convex sets and convex
functions. We write a general \emph{optimization problem} as
\[
    \begin{aligned}
        & \text{minimize} && f(x)\\
        & \text{subject to} && x \in S,
    \end{aligned}
\]
with variable $x \in \reals^n$, objective function $f: \reals^n \to \reals \cup
\{\infty\}$, and constraints $x \in S \subseteq \reals^n$. We say $x^\star$ is
a \emph{solution} to the problem if $x^\star \in S$ and $f(x) \ge f(x^\star)$
for all $x \in S$. Note that solutions may or may not exist, depending on $f$
and $S$. Alternatively, we may switch `minimize' for `maximize,' but maximizing
$f$ over $S$ is the same as minimizing $-f$ over $S$. Unlike a solution,
the \emph{optimal value} of this problem always exists and is equal to
\[
    \inf_{x \in S} f(x).
\]
For convenience, we will say that if $S \cap \dom f =\emptyset$ the optimal
value of this problem is $\infty$. (A minimizing sequence also always exists,
but we do not use this notion here.)

\paragraph{Convex problems.}
If the objective function $f$ and the set $S$ are convex, we say that the
problem is a \emph{convex optimization problem}. From before, if we also know
that the set $S \cap \dom f$ is compact and $f$ is lower semicontinuous, then
there is a solution to the problem, $x^\star \in S$. There is a very rich set
of results for problems of this form, namely the theory of duality.

\paragraph{Dual problem.} We will focus on the special case where
$S$ is an affine set and the remainder of the constraints are part of
the effective domain of $f$, which is the case we use in this paper.
This problem can be written as
\[
    \begin{aligned}
        & \text{minimize} && f(x)\\
        & \text{subject to} && Ax = b,
    \end{aligned}
\]
where $A \in \reals^{m \times n}$ and $b \in \reals^m$ are the problem data and
$x \in \reals^n$ is our problem variable. For future reference, we will call
the optimal value of this problem $s^\star$. We introduce a function called the
\emph{Lagrangian}, by relaxing the constraints to some price penalties $p
\in \reals^m$ (often called dual variables):
\[
    L(x, p) = f(x) + p^T(Ax - b).
\]
The \emph{dual function} is the optimal value of optimizing the (relaxed)
objective $L(\cdot, p)$ against fixed price penalties $p$, instead of hard
constraints:
\[
    g(p) = \inf_x L(x, p) = \inf_x \left(f(x) + p^T(Ax - b)\right).
\]
Since any point $x$ that is feasible for the original problem (\ie, satisfies
$Ax = b$) will suffer no penalty for any price $p$, we must have that
\[
    g(p) = \inf_x L(x, p) \le \inf_{Ax = b} L(x, p) = \inf_{Ax = b} f(x) = s^\star.
\]
In other words, for any price $p$, the value $g(p)$ is a lower bound to the
optimal value.

One of the most important results in convex optimization is that, under a
certain technical assumption, there exists some prices
$p^\star$ for which this inequality is tight. More specifically, if there
is a point in the relative interior of the domain of $f$, $x \in \relint \dom
f$, then there exists a set of prices $p^\star$ with
\[
    g(p^\star) = s^\star.
\]
Of course, since we know that $g(p) \le s^\star$ for all $p$, then
\[
    g(p^\star) = \sup_p g.
\]
So, to find an optimal set of prices, it suffices to find a point $p^\star$
that maximizes $g$, if we know that such prices exist. This problem:
\[
    \begin{aligned}
        & \text{maximize} && g(p)
    \end{aligned}
\]
over $p$ is known as the \emph{dual problem} (the original problem is also
referred to as the \emph{primal problem}) and has the same optimal objective
value.

\subsection{Subdifferentials}\label{app:convex-subdiff}
In general, many convex functions are not differentiable, and, unfortunately,
in many practical cases, the solution of optimization problems often lie at
points of nondifferentiability. There is a very natural extension of
differentials to convex functions that play a very similar role, called the
subdifferentials, and generalize condition~\eqref{eq:conv-ineq}.

\paragraph{Subgradients.} We say $q\in \reals^n$ is a \emph{subgradient}
of $f: \reals^n \to \reals \cup \{\infty\}$ at $x \in \reals^n$ if
\[
    f(y) \ge f(x) + q^T(y - x).
\]
for every $y \in \reals^n$. (Note the similarity between this
and~\eqref{eq:conv-ineq}.) The set of all subgradients at $x$ is called the
\emph{subdifferential} of $f$ at $x$, written $\partial f(x)$. There may be
many subgradients at some point, though in the case that $f$ is differentiable
at $x$, there will only be one. Similarly to~\eqref{eq:conv-ineq} a simple
optimality condition is, if $0 \in \partial f(x)$, then $x$ is a minimizer of
$f$, and vice versa.

In general, subgradients need not always exist everywhere, even for convex
functions. On the other hand, if $x \in \relint \dom f$, then it can be shown
that the set $\partial f(x)$ is nonempty and bounded. In almost all functions
we consider in practice, the function $f$ is indeed everywhere
subdifferentiable (has nonempty subdifferential) on its effective domain.

\paragraph{Operations on subdifferentials.}
If $f$ and $g$ are subdifferentiable at $x$, then $\partial (f+g)(x) = \partial
f(x) + \partial g(x)$. Additionally, if $f$ is defined as the supremum over a family
of functions $f_y$, written:
\[
    f(x) = \sup_{y \in Y} f_y(x),
\]
and this supremum is achieved at some $y^\star \in Y$, then it is not hard
to show that a subgradient $q$ of $f_{y^\star}$ at $x$ is a subgradient
of $f$ at $x$; \ie,
\[
    \partial f_{y^\star}(x) \subseteq \partial f(x).
\]
There are a number of other operations, but these are the only ones we will use in
this paper.

\section{Partial converse}\label{app:converse}
We will show that if $p \in \intr K$ and $\ell$ contains no line in the
direction of $p$, then there exists some $t > 0$ such that $g(tp) < g(0)$. In
other words, if there exists a strict separating hyperplane $p$ between the
sets $AX^\star$ and $Y^\star$ (which is true if, and only if, $AX^\star \cap
Y^\star = \emptyset$) then $g(0)$ is not optimal and $p$ is a
descent direction.

From before, the interior of the cone is defined
\[
    \intr K = \{p' \in \reals^m \mid p'^TAx > p'^Ty ~ \text{for all} ~ x \in X^\star, ~ y \in Y^\star \}.
\]
Now, we will show the contrapositive: if $g(tp) \ge g(0)$ for all $t > 0$, then
$p \not \in \intr K$. To see this, let $x_t \in \conv(S)$ and $y_t \in
\reals_+^n$ be the maximizers for $g(tp)$ (note that such maximizers exist
since $\conv(S)$ is a compact set and $\ell$ is lower semicontinuous and
contains no line in the direction of $p$) such that
\[
    g(tp) - g(0) = f^*(x_t) - \ell(y_t) + tp^T(y_t - Ax_t) - (\sup f^* - \inf \ell) \ge 0.
\]
Rearranging slightly, we find
\begin{equation}\label{eq:ineq-t}
    \frac{f^*(x_t) - \sup f^*}{t} + \frac{\inf \ell - \ell(y_t)}{t} + p^T(y_t - Ax_t) \ge 0.
\end{equation}
Since $f^*(x_t) \le \sup f^*$ and $\ell(y_t) \ge \inf \ell$ then
\[
    f^*(x_t) \to \sup f^* \quad \text{and} \quad \ell(y_t) \to \inf \ell,
\]
since $p^Ty_t$ and $p^TAx_t$ are bounded. (Otherwise, the left hand side of~\eqref{eq:ineq-t}
would be unbounded from below.) Additionally, inequality~\eqref{eq:ineq-t} means
\[
    p^T(y_t - Ax_t) \ge 0,
\]
for all $t$. The lower semicontinuity of $f^*$ and $\ell$ mean that $x_t \to
X^\star$ and $y_t \to Y^\star$, and, since these sequences are bounded, we must
have that there exists a subsequence $x_{t_k}$ and $y_{t_k}$ with $x_{t_k} \to
x \in X^\star$ and $y_{t_k} \to y \in Y^\star$. This, in turn, means that
\[
    p^T(y - Ax) \ge 0,
\]
so $p \not \in \intr K$.

\section{Exponential update rule}\label{app:exp-update}
In this section, we show how a log transform of the prices leads to an update
rule resembling~\eqref{eq:buterinUpdate}, proposed by Ethereum 
developers~\cite{buterin_multidimensional_2022}.
We assume that the loss function is separable and nondecreasing (\cf, \S\ref{sec:nonnegative})
and that the minimum demand condition is met (\cf, \S\ref{sec:min-demand}) so 
that the prices are strictly positive. As a result, we can write the prices $p$ 
as
\[
    p = \exp(\tilde p)
\]
for some $\tilde p \in \reals^m$, where the exponential is applied elementwise.
The gradient with respect to this new variable $\tilde p$ is 
$\nabla g(\tilde p) = \nabla g(p) \odot \exp(\tilde p)$, and the resulting 
gradient update for $\tilde p$ is 
\[
    \tilde p^{k+1} = \tilde p^k - \eta \nabla g(p^k) \odot \exp(\tilde p^k).
\]
Taking the exponential of each side, we get
\[
    p^{k+1} = p^k \odot \exp\left(-p^k \odot \eta\nabla g(p^k) \right).
\]
When we use the indicator loss function~\eqref{eq:indicator-loss}, we obtain a 
rule similar to~\eqref{eq:buterinUpdate} but with an extra factor of $p^k$ in the 
exponent:
\begin{equation}\label{eq:exp-rule}
    p^{k+1} = p^k \odot \exp(\eta (p^k \odot (Ax^k - b^\star))).
\end{equation}

\end{document}

%% file: defs.tex
\newcommand{\BEAS}{\begin{eqnarray*}}
\newcommand{\EEAS}{\end{eqnarray*}}
\newcommand{\BEA}{\begin{eqnarray}}
\newcommand{\EEA}{\end{eqnarray}}
\newcommand{\BEQ}{\begin{equation}}
\newcommand{\EEQ}{\end{equation}}
\newcommand{\BIT}{\begin{itemize}}
\newcommand{\EIT}{\end{itemize}}
\newcommand{\BNUM}{\begin{enumerate}}
\newcommand{\ENUM}{\end{enumerate}}

\newcommand{\cf}{{\it cf.}}
\newcommand{\eg}{{\it e.g.}}
\newcommand{\ie}{{\it i.e.}}

\newcommand{\ones}{\mathbf 1}

\newcommand{\reals}{\mathbb{R}}
\newcommand{\integers}{\mathbb{Z}}

\newcommand{\diag}[1]{\mathop{{\bf diag}\left({#1}\right)}}

\newcommand{\bmat}[1]{\begin{bmatrix}#1\end{bmatrix}}

\newcommand{\conv}{\mathop {\bf conv}}

\newcommand{\argmin}{\mathop{\rm argmin}}

\newcommand{\cone}{\mathop{\bf cone}}

\newcommand{\dom}{\mathop{\bf dom}}

\newcommand{\intr}{\mathop{\bf int}}
\newcommand{\relint}{\mathop{\bf rel int}}

\newcommand{\argmax}{\mathop{\rm argmax}}

\newcommand{\proj}{\mathop{\bf proj}}

\theoremstyle{definition}

\newtheoremstyle{problemstyle}  
        {3pt}                   
        {3pt}                   
        {\normalfont}           
        {}                      
        {\bfseries}             
        {\bfseries.}            
        {.5em}                  
        {\thmname{#1}\thmnumber{ #2}: \thmnote{#3}}             
\theoremstyle{problemstyle}

\newtheorem{innercustomgeneric}{\customgenericname}
\providecommand{\customgenericname}{}
\newcommand{\newcustomtheorem}[2]{%
    \newenvironment{#1}[1]
    {%
    \renewcommand\customgenericname{#2}%
    \renewcommand\theinnercustomgeneric{##1}%
    \innercustomgeneric
    }
    {\endinnercustomgeneric}
}
\newcustomtheorem{customprob}{Problem}

\newif\iftodos

\algrenewcommand\algorithmicrequire{\textbf{Input:}}
\algrenewcommand\algorithmicensure{\textbf{Output:}}

\algnewcommand\algorithmicstorage{\textbf{Storage:}}
\algnewcommand\Storage{\item[\algorithmicstorage]}

\algdef{SE}[SUBALG]{Indent}{EndIndent}{}{\algorithmicend\ }%
\algtext*{Indent}
\algtext*{EndIndent}